\documentclass[a4paper,leqno]{amsart}

\usepackage{amsmath}
\usepackage{amsthm}
\usepackage{amssymb}
\usepackage{graphicx}
\usepackage{amsfonts}
\usepackage[applemac]{inputenc}
\usepackage{mathrsfs}
\usepackage{pdfsync}

\def\C{{\mathbb C}}
\def\R{{\mathbb R}}
\def\N{{\mathbb N}}

\def\e{{\epsilon}}
\def\le{\leqslant}
\def\ge{\geqslant}

\theoremstyle{plain}
\newtheorem{theorem}{Theorem}[section]

\theoremstyle{definition}

\newtheorem{remark}[theorem]{Remark}
\newtheorem*{remark*}{Remark}

\numberwithin{equation}{section}

\begin{document}

\title[Numerical study of fractional NLS]
{Numerical study of fractional Nonlinear Schr\"odinger equations}

\author[C. Klein]{Christian Klein}
\address[C.~Klein]
{Institut de Math\'ematiques de Bourgogne
9 avenue Alain Savary, BP 47870, 21078 Dijon Cedex}
\email{christian.klein@u-bourgogne.fr}

\author[C. Sparber]{Christof Sparber}

\address[C.~Sparber]
{Department of Mathematics, Statistics, and Computer Science, M/C 249, University of Illinois at Chicago, 851 S. Morgan Street, Chicago, IL 60607, USA}
\email{sparber@math.uic.edu}

\author[P. Markowich]{Peter Markowich}
\address[P. Markowich]{King Abdullah University of Science and 
Technology (KAUST)\\ MCSE Division\\
Thuwal 23955-6900\\ Saudi Arabia}
\email{p.markowich@damtp.cam.ac.uk}

\begin{abstract}
Using a Fourier spectral method, we provide a detailed numerical investigation of dispersive Schr\"odinger type equations involving a fractional Laplacian in the one-dimensional case. 
By an appropriate choice of the dispersive exponent, both mass and energy sub- and supercritical regimes can be identified. This allows us to study 
the possibility of finite time blow-up versus global existence, the nature of the blow-up, the stability and instability of nonlinear ground states, and the 
long time dynamics of solutions. The latter is also studied in a semiclassical setting. 
Moreover, we numerically construct ground state solutions of the fractional nonlinear Schrödinger equation.
\end{abstract}

\date{\today}


\thanks{C. S acknowledges support by the NSF through grant no. DMS-1161580.  Additional support was provided through the NSF research network Ki-Net}
\maketitle

\section{Introduction}\label{sec:intro}

\subsection{Background and Motivation} This work is concerned with a 
numerical study for {\it nonlocal} dispersive equations of {\it nonlinear fractional Schr\"odinger} type (fNLS). More specifically, we consider 
equations of the form 
\begin{equation}\label{fNLS}
i \partial_t \psi = \frac{1}{2}(-\Delta)^s \psi + \gamma |\psi|^{2p} \psi, \quad \psi(0,x)=\psi_0(x),
\end{equation}
for $(t,x)\in \R\times \R^d$ and $p>0$. In addition, $\gamma = \pm 1$ distinguishes between focusing (repulsive) $\gamma =-1$ and defocusing (attractive) $\gamma =+1$
nonlinearities. Finally, the parameter $0<s\le 1$ describes the fractional dispersive nature of the equation. The 
{\it fractional Laplacian} $(-\Delta)^s$ is thereby defined via
\[
(-\Delta)^s f (x) := \mathcal F^{-1} (|k|^{2s} \mathcal F f) =  \frac{1}{(2\pi)^d} \int_{\R^d} |k|^{2s} \widehat f(k) e^{i k\cdot x} \, dk,
\]
where $\widehat f\equiv \mathcal F f $ denotes the Fourier transform of $f$. Clearly, for $s=1$ this is the usual Laplacian, whereas for $s<1$ the equation is indeed nonlocal. 
(The factor $\frac{1}{2}$ in front of the fractional Laplacian is kept for historic reasons but could be safely scaled away by replacing $x\mapsto \sqrt{2} x$)

Equation \eqref{fNLS} generalizes the classical nonlinear Schr\"odinger equation (where $s=1$), which is a canonical model for weakly nonlinear wave propagation in 
dispersive media, cf.~eg.~\cite{SS99}. In the context of quantum mechanics the case $s=\frac{1}{2}$ can be seen as a toy model for the description of particles with 
a relativistic dispersion relation $\omega(k)=\sqrt{|k|^2 +m^2}$. This has recently been used in the mathematical description of Boson-stars, see \cite{FLJ, Len}. 
Fractional NLS also arise in the continuum limit of discrete models with long range interaction \cite{KLS}, in some models of water wave dynamics \cite{IP, OS}, 
and by generalizing the Feynman path integral to include also L\'evy processes \cite{La}. 

From the mathematical point of view, fNLS equation have recently drawn quite some interest by various authors. 
For example, the question of local and/or global well-posedness of the initial value problem \eqref{fNLS} has been studied in 
\cite{GH2, GH3, GSWZ}. In addition to that, finite-time blow-up of solutions of fNLS (with Hartree type nonlinearities) has been established in \cite{CHKL, Len}. 
Moreover, the existence, uniqueness and stability properties of the associated {\it standing wave solutions} 
has been investigated in \cite{CHHO, FL, GH1}. To this end, we recall that (nontrivial) 
standing waves are obtained in the case $\gamma = -1$ by setting $\psi(t,x) = \varphi(x)e^{- i \omega t}$, $\omega \in \R$, which leads to 
the study of the following nonlinear elliptic equation:
\begin{equation}\label{stationary}
\frac{1}{2}(-\Delta)^s \varphi - |\varphi|^{2p} \varphi= \omega \varphi,
\end{equation}
see Section \ref{sec:sol} below for more details.

\subsection{Basic mathematical properties of fNLS} In this work, we are mainly interested in the interaction between the (nonlocal) dispersion and the nonlinearity in the time-evolution of \eqref{fNLS}. 
To this end, we shall take on the point 
of view that $p>0$ is fixed and $0<s\le 1$ is allowed to vary. Intuitively, we expect the model to be better behaved the stronger the dispersion, i.e. the larger $s>0$.
To obtain more insight, we first note that the following quantities are conserved by the time-evolution of \eqref{fNLS}:
\begin{equation}\label{mass}
\text{Mass:}\quad M(t) = \int_{\R^d}|\psi(t,x)|^2 \, dx =M(0),
\end{equation}
and 
\begin{equation}\label{energy}
\text{Energy:}\quad E(t) = \int_{\R^d} \frac{1}{2}|\nabla^s \psi(t,x)|^2 + \frac{\gamma}{p+1}  |\psi(t,x)|^{2p+2} \, dx =E(0) ,
\end{equation}
where $\nabla^s\psi = \mathcal F^{-1} ((-i |k|)^s \widehat \psi)$. Note that  in the defocusing (repulsive) case $\gamma =+1$, the energy is the sum of two non-negative terms 
(the kinetic and nonlinear potential energy).
This, together with the conservation of mass, allows to infer an a-priori bound on the $H^s(\R^d)$ Sobolev norm 
of $\psi$, as well as its $L^{2p+2}(\R^d)$ norm, provided that, either $\gamma = +1$ (repulsive case) or, for $\gamma = -1$ (attractive case)  the embedding $H^s(\R^d)\hookrightarrow L^{2p+2}(\R^d)$ holds. 
The latter is true for $0<p<p_*(s,d)$, where 
\begin{equation}\label{pstar}
p_*(s,d)=
\begin{cases}
\frac{2s}{d-2s}, \quad 0<s<\frac{d}{2},\\
+ \infty , \quad \, s\ge \frac{d}{2}.
\end{cases}
\end{equation}
We remark that the embedding is used for proving a local in-time existence theorem in the repulsive case.

In addition to the conservation laws above, the fNLS equation preserves the radial symmetry and is also invariant under 
the scaling transformation
\begin{equation}
\psi(t,x) \mapsto \psi_\lambda(t,x):= \lambda^{s/p}\psi( \lambda^{2s}t, \lambda x ),
    \label{scaling}
\end{equation}
for any $\lambda >0$.  In other words, if $\psi$ solves \eqref{fNLS} then so does $\psi_\lambda$. 
With this in mind, one can check that the 
under the scaling transformation \eqref{scaling}, the \emph{homogenous} $\dot H^{\sigma}(\R^d)$ Sobolev norm of $\psi_\Lambda$ behaves like
\begin{equation}
    \| \psi_\lambda \|_{\dot H^\sigma} \equiv \| \nabla^\sigma  \psi_\lambda \|_{L^2} = \lambda^{\frac{d}{2} - \sigma - \frac{s}{p}} \| \psi \|_{\dot H^\sigma} .
    \label{dotH}
\end{equation}
The equation is called $\dot H^{\sigma}$ critical whenever this scaling leaves the $\dot H^{\sigma}$ norm invariant, i.e. whenever
\begin{equation}
    \frac{d}{2} - \frac{s}{p} = \sigma .
    \label{dotHexp}
\end{equation}
For $\sigma =0$, we therefore obtain the $L^2$ critical, or {\it mass critical case} whenever the dispersion rate is $s =s^*(p,d)\equiv \frac{pd}{2}$, or, equivalently, whenever $p=\frac{2s}{d}$. 
The equation is called {\it mass subcritical} if $s>s^*$ and {\it mass supercritical} for $s<s^*$ (and vice versa for $p$). 
This should be compared to the 
situation for the usual NLS in which $s=1$ is fixed. The corresponding mass critical case is found 
for $ p=\frac{2}{d}$, particular examples being the cubic NLS in $d=2$, or the quintic NLS in $d=1$. It is well known, cf. \cite{Caz, SS99}, that in the mass subcritical case $p< \frac{2}{d}$ the 
classical NLS is globally well-posed (regardless of the sign of $\gamma$). 
On the other hand, finite time blow-up of solutions in the $\dot H^1(\R^d)$ norm can occur in the focusing case $\gamma=-1$ as soon 
as $p \ge \frac{2}{d}$.
This means that there exists a finite time $0<t^{*}<+\infty$, depending on the initial data $u_0$, such that
\[
\lim_{t\to t^{*}} \| \nabla \psi (t, \cdot)\|_{L^2} = +\infty.
\]
Moreover, it is known that for mass critical NLS, the threshold for finite time blow-up is given by the mass of the corresponding {\it ground state}, i.e. the 
unique positive radial solution $Q(x)=\varphi(|x|)$ of 
the nonlinear elliptic equation \eqref{stationary}, with $\omega =1$. 
In other words, if $p = \frac{2}{d}$ and $M(u_0)< M(Q)$, global existence still holds, whereas blow-up occurs as soon as $M(u_0)\ge M(Q)$.
For the fractional NLS an analogous dichotomy appears and has been rigorously studied in, e.g., \cite{GH2, Len}.

As we have seen there is a second conserved quantity, namely the energy. We therefore can introduce a corresponding second notion of criticality.
More precisely, the {\it energy critical case} is obtained for $s=\sigma$, in which case the kinetic energy of the solution is indeed a scale invariant quantity of the time-evolution.
This yields another critical index $s_*(p,d)=\frac{pd}{2+2p}$, which is equivalent to $p=p_*(s,d)$ as defined in \eqref{pstar}. 
Clearly, the energy critical index is always smaller than the mass critical one, i.e. $s_*<s^*$.
For classical NLS with $s=1$ fixed, the energy critical case is given by
$p=\frac{2}{d-2}$ and hence only appears in dimensions $d\ge 3$. The latter is no longer true for fractional NLS with $s<1$. In the attractive energy critical and 
supercritical case, the quantity $E(t)$ can no longer be used in order to obtain a-priori estimates on the solution. Furthermore, the classical well-posedness theory for 
semilinear dispersive PDEs breaks down as the time of existence of local solutions in general may depend on the profile $\psi_0$, not only its $H^s$ norm. 
For classical NLS ($s=1$), partial results on the existence and long time behavior of solutions in the energy critical case are 
still available, see \cite{CKGTT, DM, KV, Vi}, but a complete picture is missing so far. 
The corresponding situation for energy critical fNLS has been recently studied in \cite{GSWZ}.

\subsection{Structure of the present work} All of the above considerations paint a picture in which the theory for fNLS seems to follow closely the usual NLS results. While this is certainly true for
 basic questions such as existence and uniqueness versus finite time blow-up, the {\it nonlocal} nature of \eqref{fNLS} with $s<1$ is 
expected to have a considerable influence on more qualitative properties of the solution. 
In this paper, numerical simulations are performed in order study the influence of a nonlocal dispersion term 
on different mathematical questions, including: the particular type of finite time blow-up (e.g., self-similar or not), qualitative features of the associated 
ground states solutions (including their stability), and the possibility of 
well-posedness in the energy supercritical regime. The fact that we can vary the dispersion coefficient $0<s<1$, 
allows us to perform our simulations for both sub- and supercritical regimes in $d=1$ spatial dimension, which is a big advantage. 
(In contrast to that, numerical simulations for energy supercritical 
NLS require at least $d=3$, see \cite{CSS}.) 
For our numerical simulations we will thus fix the nonlinearity to be {\it cubic}, i.e. $p=1$, in which case \eqref{fNLS} becomes
\begin{equation}\label{fNLS1}
i \partial_t \psi = \frac{1}{2}{(-\Delta)}^s \psi + \gamma |\psi|^{2} \psi, \quad \psi(0,x)=\psi_0(x),\quad x\in \R.
\end{equation}
This equation in $d=1$ is mass critical for $s^*=\frac{1}{2}$ and energy critical for $s_*=\frac{1}{4}$. We consequently have global well-posedness 
for $\gamma =+1$ (defocusing case) and $s>\frac{1}{4}$. In the focusing case $\gamma = - 1$, we can expect finite time blow-up of solutions 
whenever $\frac{1}{2}\ge s> \frac{1}{4}$. Finally, the energy supercritical regime corresponds to $s<\frac{1}{4}$ where, in principle, a different type of blow-up
(than in the mass supercritical regime) might happen. The numerical simulations conducted are based on a {\it Fourier spectral methods}, to 
be explained in more detail in the upcoming section. Indeed, the paper is organized as follows:
\begin{itemize}
\item In Section \ref{sec:num} we will describe in more detail the numerical methods used to handle the time-evolution, as well as the steady state problem for fNLS.
\item In Section \ref{sec:Numblowup} we shall present numerical and analytical methods used for the study the blow-up phenomenon.
\item In Section \ref{sec:Qstab} we study numerically the stability of fractional ground states in the mass sub- and supercritical regime. 
\item In Section \ref{sec:blowup}, we shall study the possibility of finite time blow-up for mass critical and super-critical fNLS, by comparing the situation with the one 
occurring for the quintic and septic NLS.
\item In Section \ref{sec:energy} the time evolution of focusing and defocusing energy (super)critical fNLS numerical simulated. The possibility of blow-up is studied and 
we also study the long time behavior of the scaling invariant $\dot H^s(\R)$ norm in the defocusing case.
\item Finally, in Section \ref{sec:semi} we study the time-evolution of \eqref{fNLS1} within a semiclassical scaling regime.
\item Our main findings are summarized in Section \ref{sec:conclusio}.
\end{itemize}

\section{Numerical methods}\label{sec:num}

In the following we will discuss the numerical methods used to compute the time-evolution of the solution and its corresponding ground states.

\subsection{Numerical methods for the time-evolution}\label{sec:numtime}

For the numerical integration of \eqref{fNLS}, we use
a Fourier spectral method in $x$. The reason for this choice is that 
the fractional derivatives are most naturally computed in frequency 
space which is approximated via a discrete Fourier transform computed 
through an FFT (fast Fourier transform). The excellent approximation properties 
of a Fourier spectral method for smooth functions are also extremely 
useful. This is especially important in the context of dispersive 
equations since spectral methods are known for a minimal  numerical 
dissipation which (in principle) could overwhelm dispersive effects 
within our model. 

\begin{remark} For the same reasons, a recent numerical study using the same numerical methods 
has been conducted for fractional KdV and BBM type equations, see \cite{KS14}.\end{remark}

The discretization in Fourier space leads to a system of 
(stiff) ordinary differential equations for the Fourier 
coefficients of $\psi$ of the form
\begin{equation}
   \partial_t \widehat{\psi}=\mathcal{L}\widehat{\psi}+\mathcal{N}(\psi),
    \label{fnlsfourier}
\end{equation}
where $\mathcal{L}=-i|k|^{2s}/2$ and where
$\mathcal{N}(\psi)=i\gamma \widehat{|\psi|^{2p}\psi}$ 
denotes the nonlinearity. It is an advantage of Fourier methods 
that the $x$-derivatives and thus the operator $\mathcal{L}$ are 
diagonal. For equations of 
the form (\ref{fnlsfourier}) with diagonal $\mathcal{L}$, there are many efficient high-order time 
integrators. In particular, the performance of several fourth order methods was recently 
compared in \cite{etna} by using the (semiclassically scaled) cubic NLS as a bench mark. 
It was shown that in the defocusing case, a 
\emph{time-splitting scheme} as in \cite{BJM,BJM2} was the most efficient, 
whereas in the focusing case a {\it composite Runge-Kutta method} 
\cite{Dris} is preferable. The two codes are also used to test 
each other and were found to agree within the indicated numerical 
precision. We shall therefore use these two approaches also in this 
paper.

\begin{figure}[htb!]
   \includegraphics[width=0.49\textwidth]{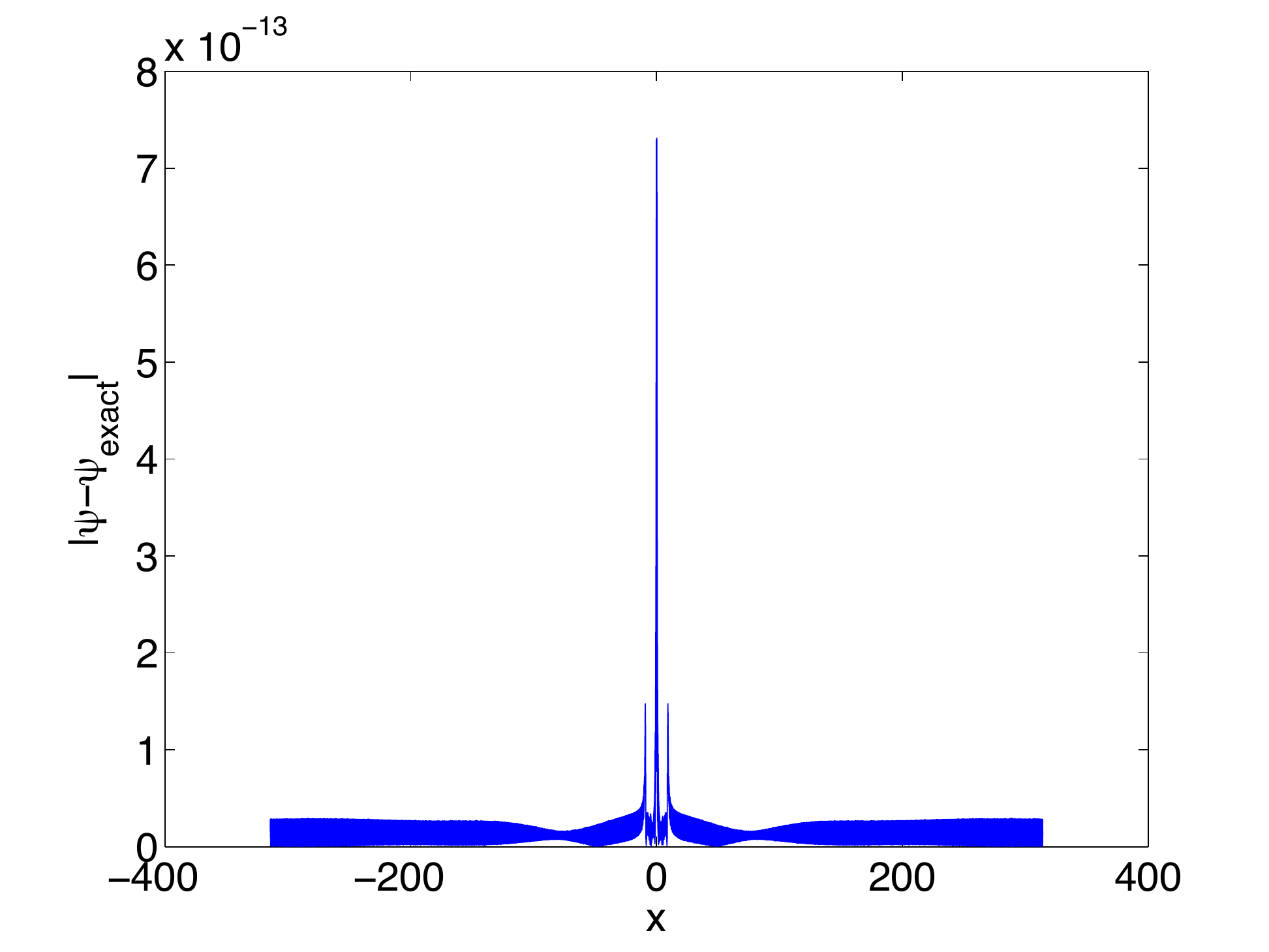}
   \includegraphics[width=0.49\textwidth]{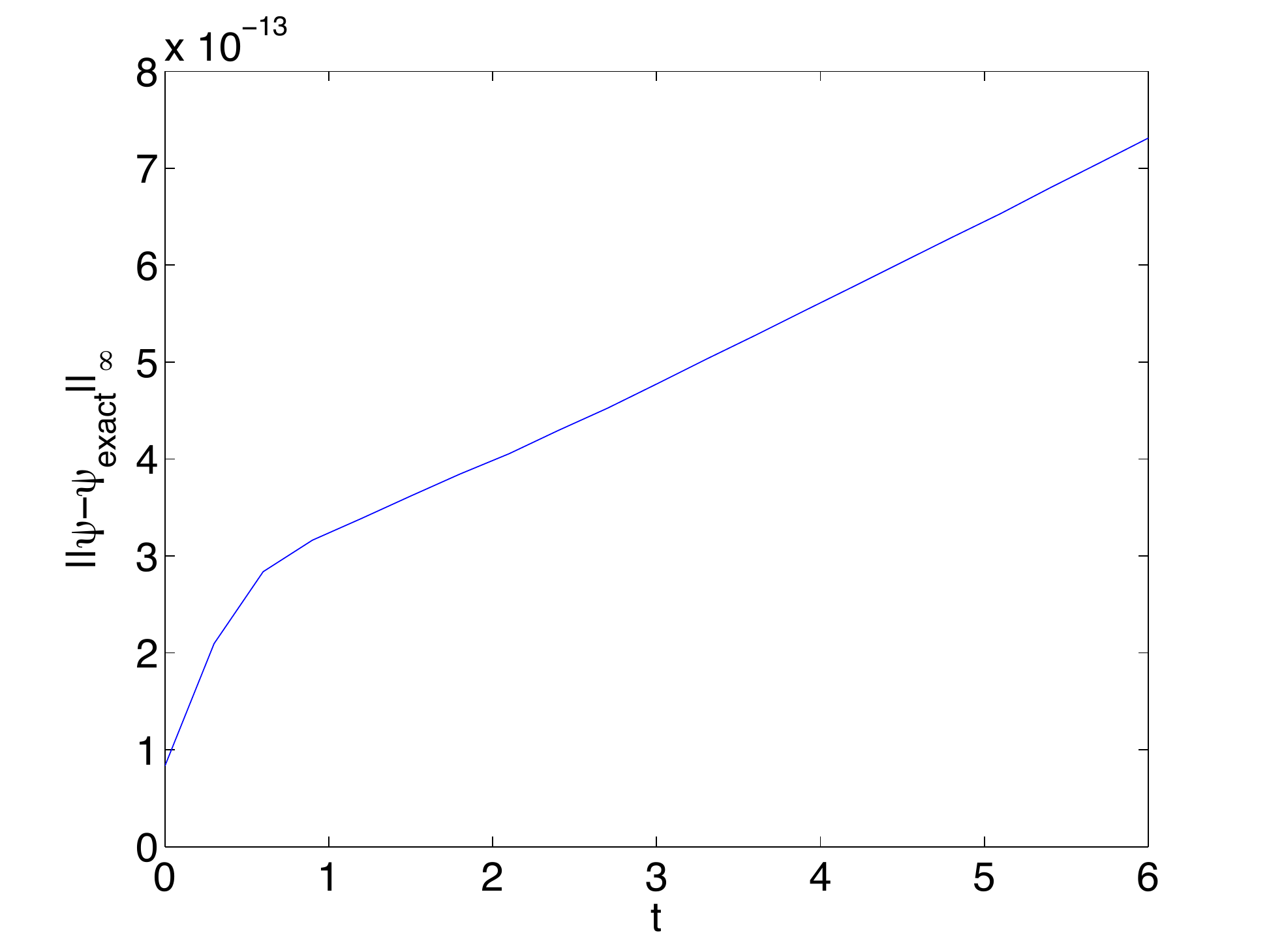}
 \caption{Left: Difference between the numerical solution of the fNLS 
 equation (\ref{fNLS1}) with $s=0.6$ and $\psi_{0}=Q$ and the 
 exact solution $Q(x)e^{it}$, at $t=6$. Right:
 $L^{\infty}$ norm of the difference in dependence of time.}
 \label{nlsfracsoldiff}
\end{figure}

In order to test the numerical methods, we take $\psi_0(x)=Q(x)$, i.e. the ground state solution whose 
numerical construction is explained in the next subsection. 
In this case, the time-dependence of the 
exact solution of \eqref{fNLS} is simply given by $\psi(t,x)=Q(x)e^{it}$. A comparison of the 
numerical solution of the fNLS with initial data $\psi_{0}=Q$ therefore tests 
both $Q$ and the time-evolution. 
In Fig.~\ref{nlsfracsoldiff} we take $p=1$ (cubic nonlinearity), 
$s=0.6$ and show the difference between the numerical solution and $Q(x)e^{it}$ for 
$N=2^{16}$ and $N_{t}=20000$ time steps for $t\leq 6$. It can be seen 
that the ground state is reproduced up to errors of order $10^{-12}$, i.e., 
essentially with machine precision (which is roughly $10^{-14}$ in our case).

In our numerical computations, we also ensure that 
the computed relative energy of the solution, i.e.
\begin{equation}
    \Delta_E=\Big| \frac{E(t)}{E(0)}-1\Big|
    \label{Delta}
\end{equation}
remains small up to the threshold $\Delta_E <10^{-3}$. For the example 
in Fig.~\ref{nlsfracsoldiff}, the quantity $\Delta_E$ is of the order 
of machine precision in accordance with expectations. Generally this 
quantity is smaller than $10^{-10}$ in our computations unless 
otherwise noted.

\subsection{Numerical construction of fractional ground states} \label{sec:sol}

Recall that standing wave solutions of \eqref{fNLS} are obtained in the focusing case $\gamma=-1$ 
by setting $\psi(t,x) = \varphi(x)e^{- i \omega t}$ for some $\omega \in \R$. 
By rescaling
\begin{equation}
 \varphi_{\omega}=\omega^{1/(2p)} \varphi_1(x\omega^{1/(2s)}),
    \label{Qresc}
\end{equation}
we can w.r.o.g. assume $\omega=1$ and hence $\varphi\equiv \varphi_1$ solves 
\begin{equation}
\frac{1}{2}(-\Delta)^{s}\varphi+\varphi=|\varphi|^{2p}\varphi
    \label{soliton0}.
\end{equation}
Solutions $\varphi \in H^s(\R^d)\cap L^{2p+2}(\R^d)$ of this equation exist for admissible $0<p<p^*$, where $p_*=p_*(s,d)$ is given by \eqref{pstar}.
Indeed, by invoking  Pohozaev-type identities, it can be shown that equation \eqref{soliton0} does not admit any nontrivial solution in $H^s(\R^d)\cap L^{2p+2}(\R^d)$, when $p\ge p_*$.
Of special interest are solutions with minimal energy, so-called \emph{ground 
states}, which are known to be real and radially symmetric and thus satisfy
\begin{equation}
\frac{1}{2}(-\Delta)^{s}Q+Q=Q^{2p+1}
    \label{soliton}.
\end{equation}
For $s\not =1$, ground state solutions $Q$
decay like $|x|^{-(d+2s)}$ as $|x|\to\infty$, i.e., only algebraically fast, see \cite{FL}. This is in contrast to the case of the usual NLS ground states 
obtained for $s=1$, which are known to decay exponentially fast, see, e.g., \cite{Caz}. Indeed, for $s=d=1$, one has the well-known explicit 
solution of the so-called bright solitary wave (at time $t=0$):
\begin{equation}\label{Q}
Q(x) = \left( \frac{p+1}{\cosh^2(\sqrt{2} p x)}\right)^{1/(2p)}, \quad 0<p<p^*.
\end{equation}
For fNLS, with $s<1$, no explicit solutions of \eqref{soliton} are known.

To solve equation (\ref{soliton}) numerically, we use the same 
approach as in \cite{KS14} to which we refer the reader for more details. 
The basic idea is to expand $Q$ on a finite interval $x\in D[-\pi,\pi]$, $D>0$, in 
a discrete Fourier series, computed via FFT. 
 In Fourier space, equation (\ref{soliton}) takes the form
\begin{equation}
 F(Q):=   \left(\frac{1}{2}|k|^{2s}+1\right)\widehat{Q}-\widehat{Q^{2p+1}}=0
    \label{solfourier},
\end{equation}
subject to periodic boundary conditions.
Due to the slow algebraic decay of $Q$ for $s<1$, the constant $D$ thereby has to be chosen sufficiently 
large in order to reduce the discontinuity of the derivatives of $Q$ at the boundaries 
of the computational domain. It is well known that such 
discontinuities imply an algebraic decrease of the Fourier 
coefficients with the wave number $k$ and thus a slow convergence of 
the numerical approximation with the number $N$ of Fourier modes. We choose here 
$D=100$ and $N=2^{16}$ Fourier modes. Numerically, $\widehat{Q}$ is approximated by a discrete Fourier 
transform, i.e., by a finite vector, which implies that a system of 
$N\gg 1$ nonlinear equations has to be solved. To this end, we invoke an iterative Newton method in order to find 
the (nontrivial) zeroes of the function $F(Q)$. This means that we iterate
\[
\widehat Q_{n+1} = \widehat Q_n - J^{-1}(F)(\widehat Q_n),
\]
where $J$ is the Jacobian of $F$.
However, for $N=2^{16}$, the dimension of the Jacobian is too high 
to be efficiently implemented, and we therefore apply a 
Newton-Krylov method. This means that the inverse of the Jacobian is 
computed via GMRES \cite{gmres}, iteratively. By doing so only the action of the 
Jacobian on a vector has to be computed, whereas the full Jacobian is 
never explicitly stored. 

\begin{figure}[htb!]
   \includegraphics[width=0.7\textwidth]{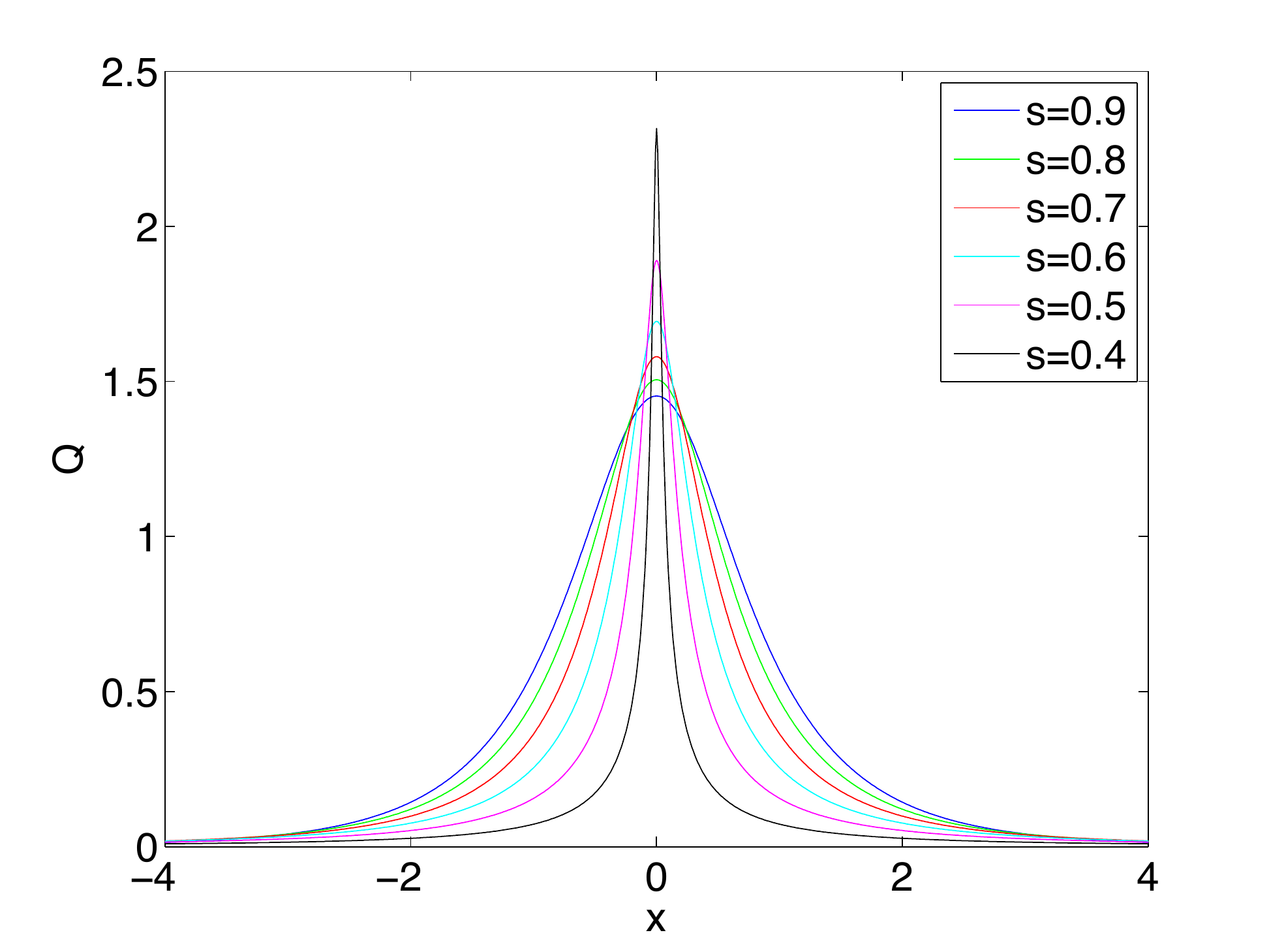}
 \caption{Ground state solutions $Q(x)$ of \eqref{fNLS1} for different values of $s$, depicted close to the origin.}
 \label{solfig}
\end{figure}

An additional obstruction is given by the fact that \eqref{soliton}, or equivalently \eqref{solfourier}, always has the 
trivial solution $\widehat{Q}=0$. Thus a fixed point iteration in general will converge to the latter. To circumvent this 
problem, we have to make sure to start sufficiently close to the exact nontrivial solution $Q$.
For $s=1$, the latter is given explicitly by \eqref{Q} and thus, we 
use a continuation method, i.e., we start with values of $s$ close to $1$, say $s=0.9$,
and an initial iterate given by \eqref{Q}. Then we 
use the solution for $s=0.9$ as the starting point for the iteration for $s=0.8$ 
and so on. The results can be seen in Fig.~\ref{solfig} where we have chosen $p=1$, i.e. a cubic nonlinearity. 
We see that the smaller $s<1$, the more the ground state solution becomes peaked, and 
the slower its spatial decay as $|x|\to\infty$ in accordance with 
the theoretical predictions.

\begin{figure}[htb!]
   \includegraphics[width=0.49\textwidth]{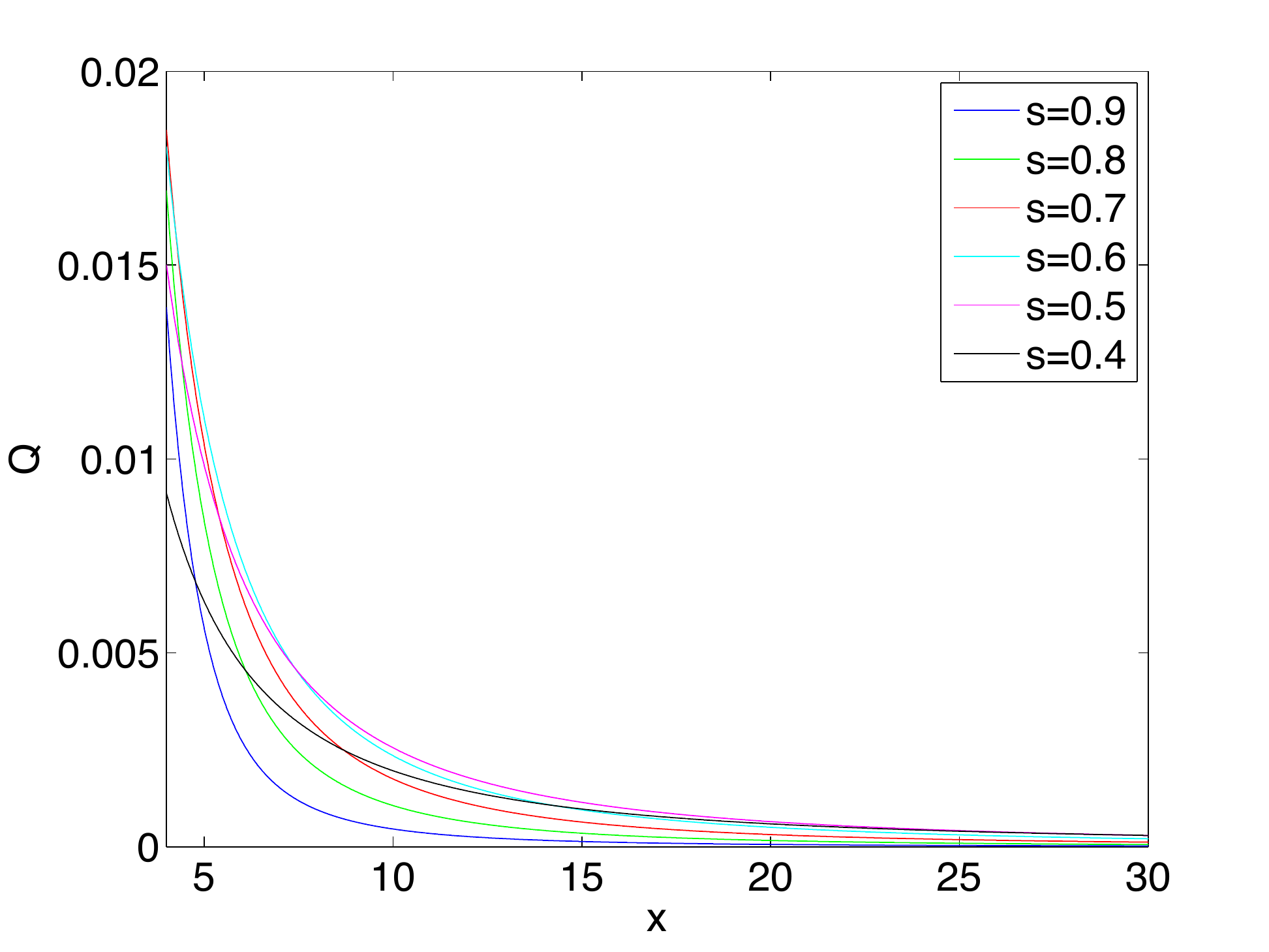}
   \includegraphics[width=0.49\textwidth]{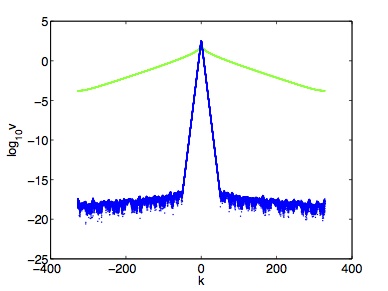}
 \caption{Left: The asymptotic behavior for large $x$ of the ground state solutions $Q$ for different values of $s$. 
 Right: The modulus of the Fourier coefficients for two of the ground 
 states shown in Fig.~\ref{solfig}, in blue for $s=0.9$, in 
 green for $s=0.4$.}
 \label{nlsfracsolfourier}
\end{figure}

\begin{remark}The slow decay of $Q$ also affects the convergence of the 
iteration, and hence we have to decrease $s<1$ in smaller and smaller increments in order to to assure convergence. 
In each case, the iteration is stopped whenever equation \eqref{solfourier} is satisfied to better than $10^{-12}$.
This implies that the solutions are 
well resolved in Fourier space for larger $s>0.5$. It can be seen in Fig.~\ref{nlsfracsolfourier} that the 
modulus of the Fourier coefficients decreases to machine precision 
for the high wave numbers, whereas they only decrease to $10^{-4}$ 
for $s=0.4$. 
\end{remark}

\begin{figure}[htb!]
   \includegraphics[width=0.49\textwidth]{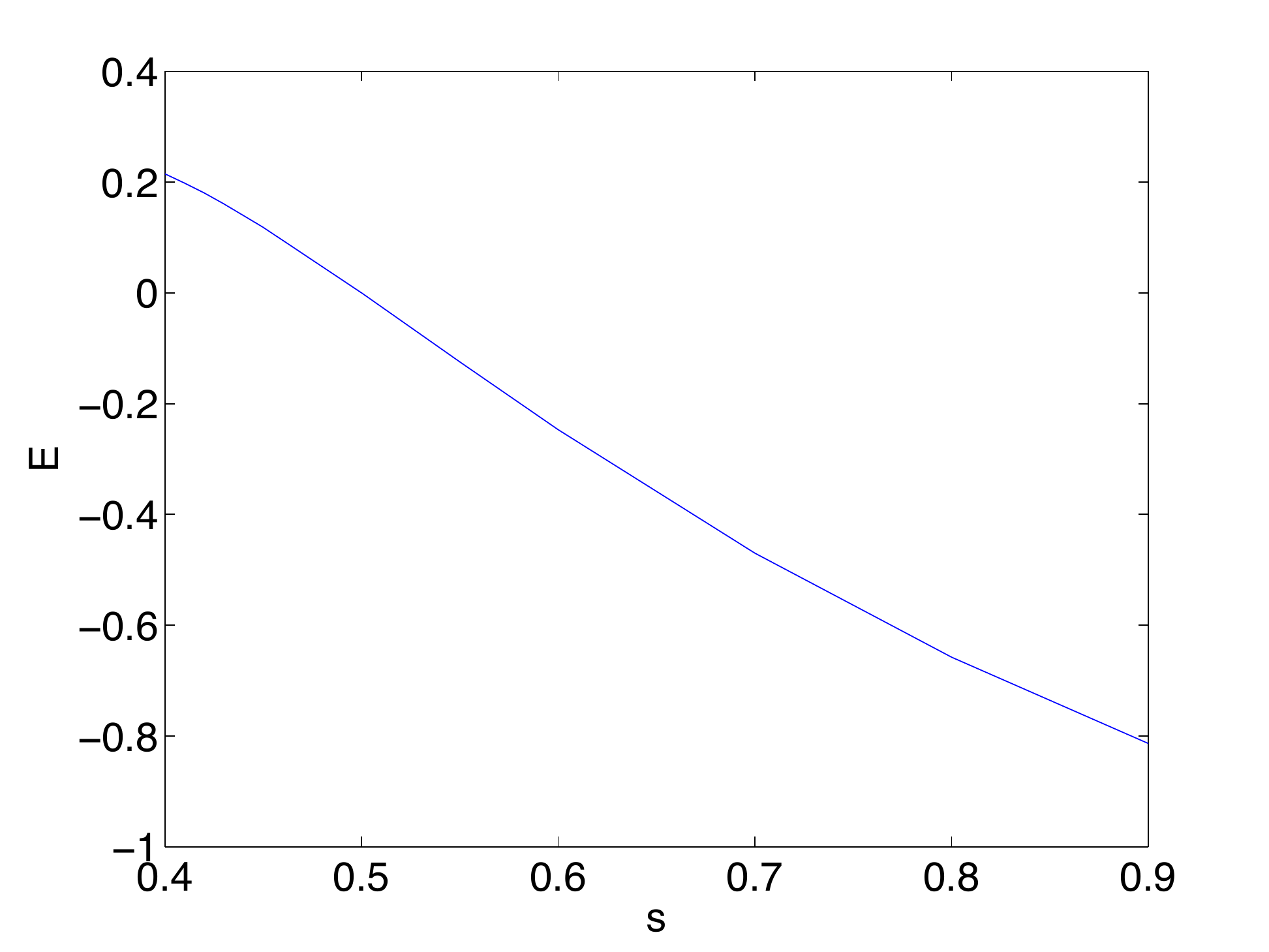}
   \includegraphics[width=0.49\textwidth]{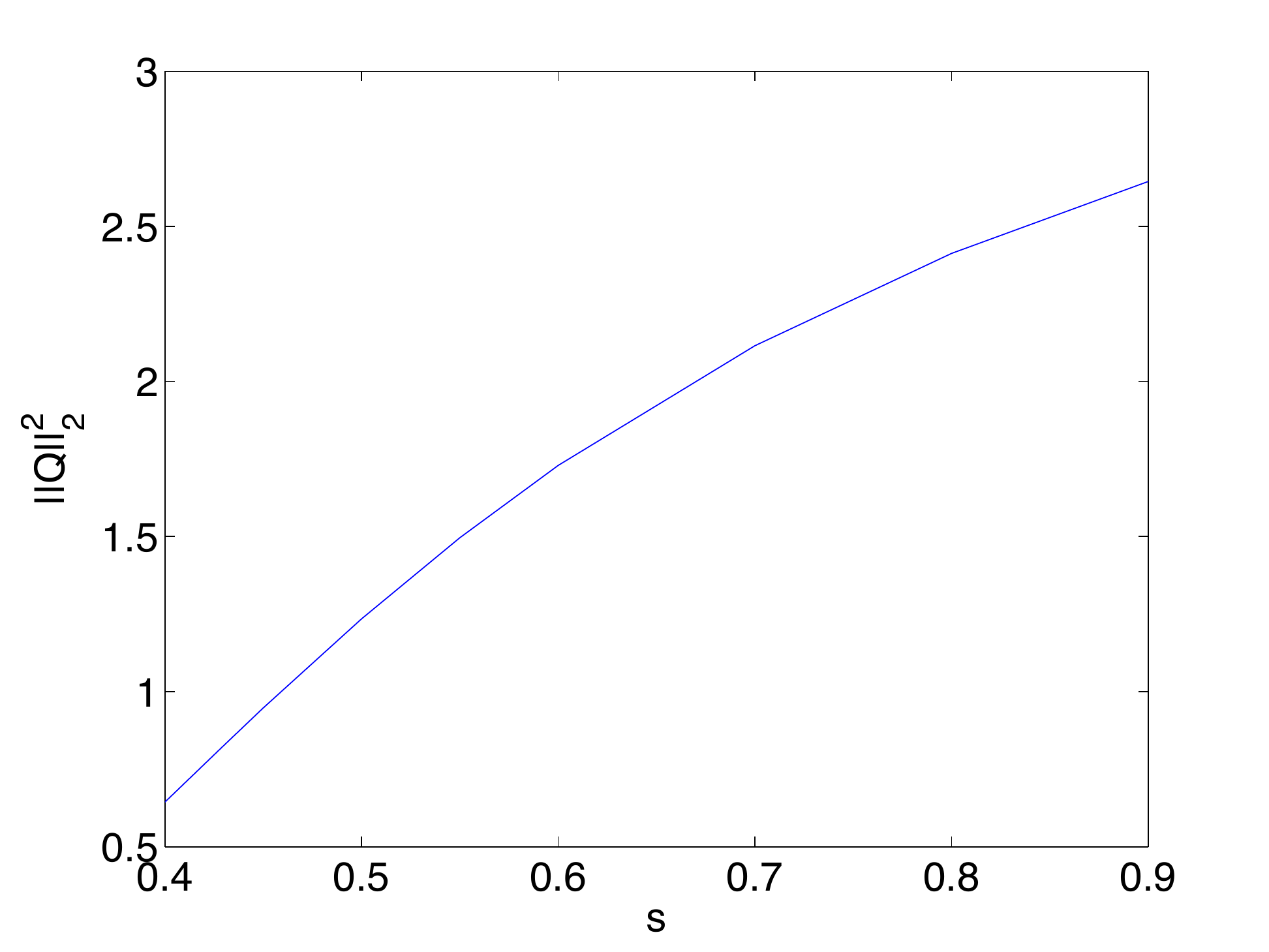}
 \caption{The $s$ dependence of the total energy (left) and mass (right) of the 
 ground state solutions $Q(x)$ to \eqref{fNLS1}.}
 \label{nlsfracsolEM}
\end{figure}

Finally, the total energy and the mass of the ground solutions $Q$ are depicted in dependence of $s$ in 
Fig.~\ref{nlsfracsolEM}. It can be seen that the energy is 
monotonically decreasing with $s$ while the mass is increasing. 
Note that the energy vanishes with numerical precision 
($\approx 10^{-5}$) in the $L^{2}$ critical case, where $p=1$ and $s=0.5$.


\section{Methods for the numerical study of blow-up}\label{sec:Numblowup}
In this section we briefly present the methods used for the numerical 
study of finite-time blow-up. Firstly we use a dynamic rescaling which allows in 
principle an adaptive mesh refinement near blow-up. Secondly, we 
explain how to numerically identify singularities on the real axis by 
tracing singularities of the solution in the complex plane via the 
asymptotic behavior of the Fourier coefficients. 

\subsection{Dynamical rescaling}
In the numerical study of blow-up in NLS equations, dynamically 
rescaled codes have proven to provide an interesting approach, see 
\cite[Chapter 6]{SS99} and the references therein. In the case of {\it radially symmetric} solutions $\psi(t,r)\equiv \psi(t,|x|)$ 
one thereby introduces the change of variables 
\begin{equation}
    y = \frac{r}{L(t)},\quad 
    \frac{d\tau}{dt}=\frac{1}{L^{2s}(t)},\quad \psi(t,r) = L(t)^{-s/p}\Psi(\tau, y)
        \label{resc}.
\end{equation}
Using this for the focusing ($\gamma = -1$) fNLS, we find a rescaled equation for $\Psi$ in the form
\begin{equation}
    i\partial_\tau \Psi =i a(\tau) \left(\frac{s}{p}\Psi+y \partial_y \Psi\right)
    +\frac{1}{2}(-\Delta)^{s}\Psi-|\Psi|^{2p}\Psi
    \label{fnlsres}.
\end{equation}
where $(-\Delta)^{s} \Psi (x):= \mathcal F^{-1}({|k|^{2s}\widehat \Psi(k)})$, for $k\in \R$, and
\begin{equation}\label{a}
a= L^{2s-1}\frac{dL}{dt}= \frac{d \ln  L} {d\tau}.
\end{equation}
Under this rescaling the mass and energy behave like
\begin{equation}
    M=  
    \frac{|\mathbb S^{d-1}|}{L(\tau)^{d-2s/p}}\int_0^\infty |\Psi(\tau, y)|^{2} \,y^{d-1} \, dy,
    \label{Mresc}
\end{equation}
and
\begin{equation}
    E=  
    \frac{|\mathbb S^{d-1}|}{L(\tau)^{d-2s-2s/p}}\int_0^\infty \left(\frac{1}{2} |\partial_y^{s}\Psi(\tau, y)|^{2}
   -\frac{1}{p+1}|\Psi(\tau, y)|^{2p+2}\right) y^{d-1}\, dy,
    \label{Eresc}
\end{equation}
where $\partial_y^s f (x):= \mathcal F^{-1}({(-ik)^s\widehat f(k)})$.
The scaling function $L$ should be chosen in such a way that $L(\tau)\to 0$, sufficiently fast, as $\tau \to +\infty$. 
It is then expected that, as $\tau \to +\infty$, both $a\to a^\infty$ and $\Psi \to \Psi^\infty$ become $\tau$-independent (in the mass supercritical case).  
The profile $\Psi^\infty$, consequently solves
\begin{equation}
   i a^{\infty}\left(\frac{s}{p} \Psi^{\infty}+y \partial_y \Psi^{\infty}\right)
    +(-\Delta_{y})^{s}\Psi^{\infty}-
   |\Psi|^{2p}\Psi=0, \quad y\in \R_+.
    \label{fnlsrescinfty}
\end{equation}
In $d=1$, this is a fractional ordinary differential equation. It is 
not known whether this equation has localized solutions, and if yes, 
whether these are unique. If such a unique solution exists, it will 
give the blow-up profile of the self similar blow-up. It is beyond 
the scope of the present article to address the related problems. 
\begin{remark}
If $a^{\infty}=0$, this equation corresponds to the standing wave equation \eqref{stationary} with $\omega=0$. The associated 
stationary solution $\varphi=W(x)$, often called Rubin-Talenti solution, plays a similar role for energy critical NLS, as does the ground state solution 
$\varphi=Q(x)$ for mass critical NLS, see \cite{DM} for more details.
\end{remark}
There are different ways of constructing $L(t)$, see \cite{SS99}. One 
of them invokes the use of  an integral norm of $\psi$ which goes to infinity at the blow-up. 
This is  preferable from the numerical point of view and hence, we shall choose a scaling which keeps $\| \partial_y \Psi(\tau, \cdot)\|_{L^2}$ constant. 
This leads to 
\begin{equation}\label{L}
L(t)^{1-d/2 + s/p} = \left(\frac{ \| \partial_y \Psi_0\|_{L^2}}{ \| \partial_r \psi(t, \cdot)\|_{L^2}}\right)
\end{equation}
where the constant $\| \partial_y \Psi_0\|_{L^2}$ is chosen to be 
$ \| \partial_r \psi_0\|_{L^2}$. For given  $\| \partial_r \psi(t, 
\cdot)\|_{L^2}$, we can read off  the time dependence of $L$ from 
(\ref{L}). ALternatively to obtain an equation for $a$ from 
(\ref{L}), we differentiate $ \| \partial_y \Psi_0\|_{L^2}$ knowing 
the resulting expression vanishes by assumption. Then we use 
(\ref{fnlsres}) to eliminate the $\tau$-derivative of $\Psi$ which 
leads to an equation involving $a$. After some partial integrations, 
we end up with 
\begin{equation}
    a(\tau)=\frac{2|\mathbb S^{d-1}|}{(2s/p+1)\|\partial _{y} \Psi \|_2^{2} }\int_0^\infty |\Psi|^{2p} \, \text{Im} \, (\bar{\Psi}\partial^2_y\Psi ) \, y^{d-1}\,dy
    \label{L2x}.
\end{equation}
This allows us in principle to study the type of the 
blow-up for fNLS in a similar way as it has been done for generalized Korteweg-de 
Vries equations in \cite{KP2013}. 
But it was shown 
numerically in \cite{KP2013} that generic rapidly decreasing hump-like 
initial data lead to a 
tail of dispersive oscillations as $|x|\to\infty$ with slowly 
decreasing amplitude. Due to the imposed periodicity (in our numerical domain), these 
oscillations reappear after some time
on the opposing side of the computational domain 
and lead to numerical instabilities in the dynamically rescaled 
equation. The source of these problems is the term $y\Psi_{y}$ in 
(\ref{fnlsres}) since $y$ is large at the boundaries of the 
computational domain. Therefore this term is very sensitive 
to numerical errors.  For gKdV this can be addressed 
by using high resolution in time and large computational domains. It 
turns out that for fractional KdV equations, see 
\cite{KS14}, and for fNLS, the dispersive oscillations have an 
amplitude that decreases very slowly towards infinity which is also 
reflected by the slow decrease of the solitons. The consequence of 
this is that  we 
cannot compute long enough with the dynamically rescaled code 
to get conclusive results. Instead we 
integrate fNLS directly, as described above, and then we use  
post-processing to characterize the type of blow-up via the above 
rescaling. For instance, we read off the 
time evolution of the quantity $L$ from (\ref{L}).

Under the hypothesis that $L(\tau) \sim\exp(-\kappa \tau)$ with $\kappa>0$ some positive constant, (\ref{resc}) yields a connection between $t$ and $\tau$. 
Namely, 
\begin{equation}
    L(t) \propto (t^{*}-t)^{\frac{1}{2s}}
    \label{Lt},
\end{equation} 
where $t^*>0$ is the blow-up time, corresponding to $\tau = +\infty$.
With (\ref{L}) and (\ref{resc}), this 
implies
\begin{equation}
    \| \partial_{x} \psi(t,\cdot) \|_{2}^{2}\propto (t^{*}-t)^{-(1/p+1/(2s))},\quad
    \|\psi(t,\cdot) \|_{\infty}\propto 
    (t^{*}-t)^{-1/(2p)}.    
    \label{genscal}
\end{equation}
In particular, for $s=1$ we have $L\propto \sqrt{t^{*}-t}$, which is the expected blow-up rate for NLS in the mass supercritical regime.
In the mass critical case $p=2$, one finds a correction to \eqref{L} in the form 
\begin{equation}
L(t) \propto \sqrt{\frac{t^{*}-t}{\ln |\ln (t^{*}-t)|}},
    \label{L2scal}
\end{equation}
i.e., one has $\tau\propto \ln (t^{*}-t)(1-\ln |\ln (t^{*}-t)|)$ 
instead of $\tau\propto \ln (t^{*}-t)$. This so-called log-log-scaling regime for mass critical NLS has been rigorously proved in \cite{MR}.
We will test whether such scalings can be observed in the 
numerical experiments for fNLS, but it cannot be expected that the 
logarithmic corrections can be seen numerically.

\subsection{Singularity tracing in the complex plane}

In the case of a finite-time blow-up, we observe essentially two types of 
behavior of the numerical solution. Either the $L^{\infty}$ norm of 
the solution becomes so large that the computation of the nonlinear 
terms in the fNLS equation leads to an overflow error. In this case 
the code breaks down by producing $NaN$ results. The other possibility is that 
the code runs out of resolution in Fourier space which is indicated by a deterioration of 
the Fourier coefficients. The latter allows for an identification of 
an appearing singularity as follows (see also \cite{KR2013a, KR2013b, 
SSF}): The function $\psi$ to be studied on the real axis is assumed to have 
a continuation in the complex plane in a strip around the 
real axis denoted by $f(z)$.
We recall the fact that, in the complex plane, a (single) essential singularity $z_0\in \C$ of a 
function $f$, such that $f(z) \sim (z-z_{0})^{\mu}$ for $z\sim 
z_{0}$, with $\mu\not \in \mathbb Z$, results in the following 
asymptotic behavior (for $|k|\to\infty$) for the corresponding 
Fourier transform (see e.g.~\cite{asymbook})
\begin{equation}
    |\widehat{f}(k)|\sim 
    \frac{1}{k^{\mu+1}} e^{-k\delta},\quad |k| \gg 1,
    \label{fourasymp}
\end{equation}
where $\delta=\text{Im}\, z_{0}$. The quantity $\mu$ thereby 
characterizes the type of the singularity.

In \cite{KR2013a, KR2013b}  this 
approach was used to quantitatively identify the time where the 
singularity hits the real axis, i.e., where the real solution becomes 
singular, since it was shown that the quantity $\delta$ can be reliably identified 
from a fitting of the Fourier coefficients. Unfortunately, this is not true for $\mu$, for which the numerical inaccuracy is too large. 
In the case of focusing NLS, it was shown in \cite{KR2013b} that 
the best results are obtained when the code is stopped once the 
singularity is closer to the real axis than the minimal resolved 
distance via Fourier methods, i.e.,
\begin{equation}
    m:=2\pi \frac{D}{N},
    \label{mres}
\end{equation}
with $N\in \N$ being the number of 
 Fourier modes and $2\pi D$ the length of the computational domain in 
 physical space.  All values of $\delta <m$ 
 cannot be distinguished numerically from $0$. 
 
 Note that the 
 time at which the code is stopped because of the criterion above is not 
 the same as the blow-up time itself. Rather, it is only the time where the code stops to be 
 reliable. As mentioned above, we will always provide sufficient resolution in 
 time so that that only the lack of resolution in Fourier space makes the code 
 stop. The blow-up time will then be determined from the numerical 
 data by fitting to the scalings given in the previous subsection. We generally choose the time step in blow-up scenarios such that the 
accuracy is limited by the resolution in Fourier space, i.e., that a 
further reduction of the time step for a given number of Fourier 
modes does not change the final result within numerical accuracy.

\section{Stability of ground states}\label{sec:Qstab}

In this section we will study perturbations of the ground state solutions 
constructed before. This is done for cubic nonlinearities $p=1$ and different 
values of the parameter $s$. To this end, we will consider initial 
data for \eqref{fNLS1} 
of the form
\begin{equation}\label{iniQ}
\psi_0(x) = \alpha Q(x), \quad \alpha \in \R,
\end{equation}
where $Q$ is the ground state solution determined numerically as described in 
Section~\ref{sec:sol}. The factor $\alpha$ will be either chosen to 
be a constant $\alpha \approx 1$ or to be an $x$ dependent phase. 
Note that qualitatively similar results as shown here are also found for 
localized perturbations of the form: $\psi_0(x) =  Q(x) +\varepsilon \, e^{-|x|^2}$, with $\varepsilon <1$. 

\subsection{Perturbed ground states in the mass subcritical regime}

It is known \cite{CHHO} that the ground state solutions are stable in the mass 
subcritical case, i.e. $s>\frac{1}{2}$. In fact, if we propagate initial data of the form \eqref{iniQ} we find that the perturbed ground state starts to 
oscillate around what appears to be a stationary solution with frequency $\omega\in \R$, i.e. we find that $\psi(t,x) =  Q_{\omega}(x)e^{i\omega t}$. This can 
be seen in Fig.~\ref{sol_{09}_{09}}, where we have solved the initial 
value problem \eqref{fNLS1} subject to data \eqref{iniQ} with $\alpha=0.9$. 
We use $N=2^{16}$ Fourier modes for $x\in100[-\pi,\pi]$ 
and $N_{t}=10^{4}$ time steps for $t<30$. It can be seen that the 
initial hump decreases in height and then starts to exhibit damped 
oscillations around what appears to be a rescaled ground state function. This 
is reminiscent of the so-called {\it breather solutions} known for classical NLS.
\begin{figure}[htb!]
   \includegraphics[width=0.7\textwidth]{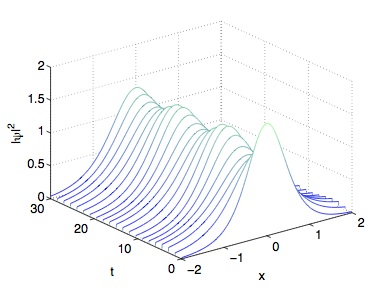}
 \caption{Modulus squared of the solution to the focusing fNLS equation 
 (\ref{fNLS1}) with $s=0.9$ for initial data $\psi_{0}(x)=0.9Q(x)$.}
 \label{sol_{09}_{09}}
\end{figure}

The damped oscillations around some presumably constant 
asymptotically constant are clearly visible if one looks at the 
$L^{\infty}$ norm of the solution, see Fig.~\ref{sol_{09}_{09}max}. 
Since the $L^{2}$ norm of the solution is a conserved quantity, the scaling (\ref{scaling}) allows us to infer 
a bound on $\omega$, given by
$\omega^{1-1/(2s)}\leq\alpha^{2}$. For $\alpha=0.9$ this would imply 
that the $L^{\infty}$ norm of the ground state with the maximal 
$\omega$ would be roughly equal to $1.146$. Fig.~\ref{sol_{09}_{09}max} suggests 
that this is indeed the amplitude of the final state. This would mean that 
the ground state is stable, and that a perturbed ground state leads 
asymptotically for large $t$ to a steady state with the mass of the 
perturbed state. In the same figure we show the $L^{\infty}$ norm of 
the solution for the fNLS equation for initial data \eqref{iniQ} with $\alpha =1.1$. There 
are much more oscillations in this case, but the final state appears 
to have an $L^{\infty}$ norm of roughly 1.8 (the maximal possible 
$L^{\infty}$ norm of the ground state having the same mass as the 
initial data would be $\approx 1.800$). Thus also in this case the 
final state appears to be a stationary solution corresponding to the mass 
of the initial data. 
\begin{figure}[htb!]
   \includegraphics[width=0.49\textwidth]{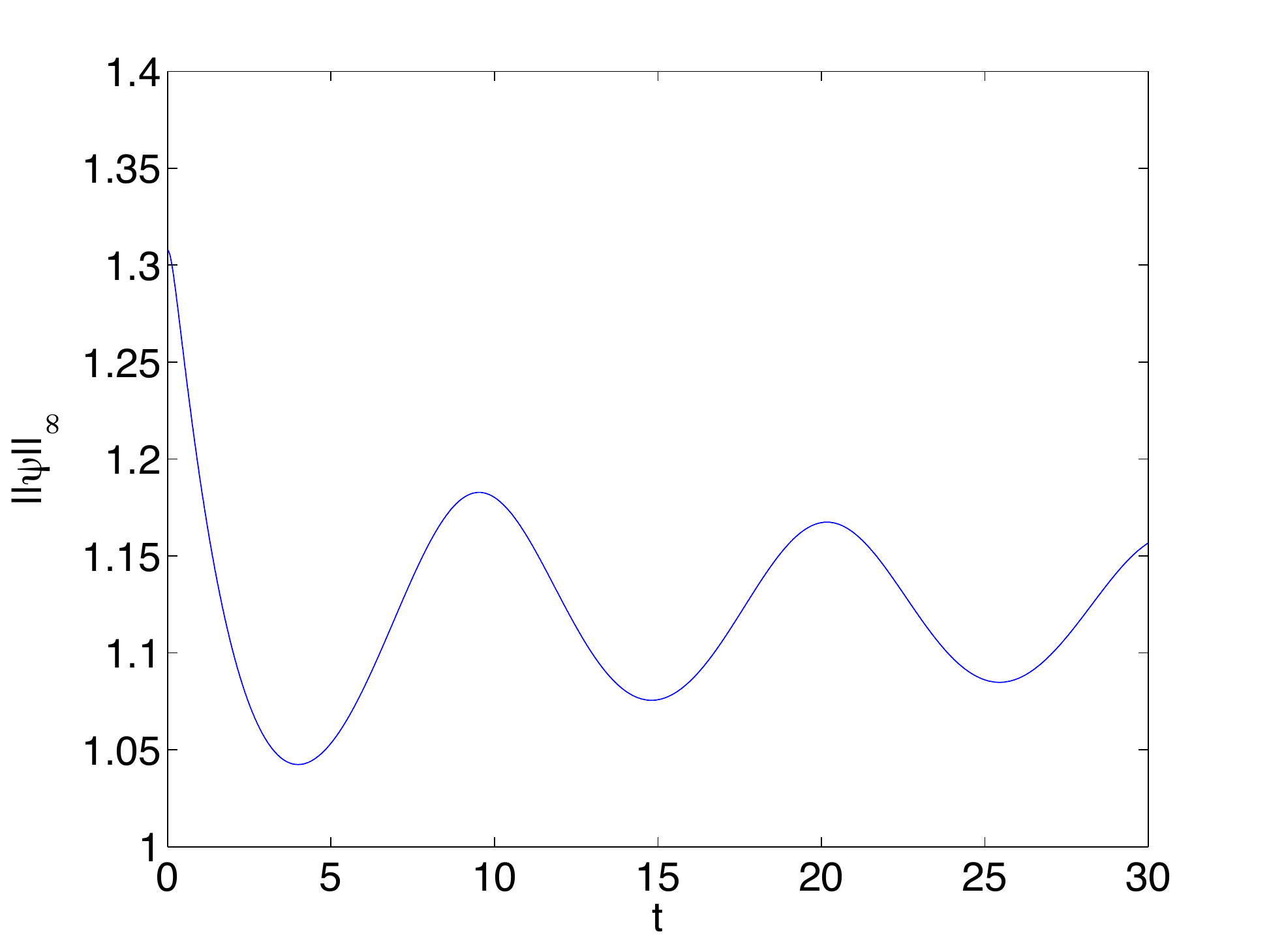}
   \includegraphics[width=0.49\textwidth]{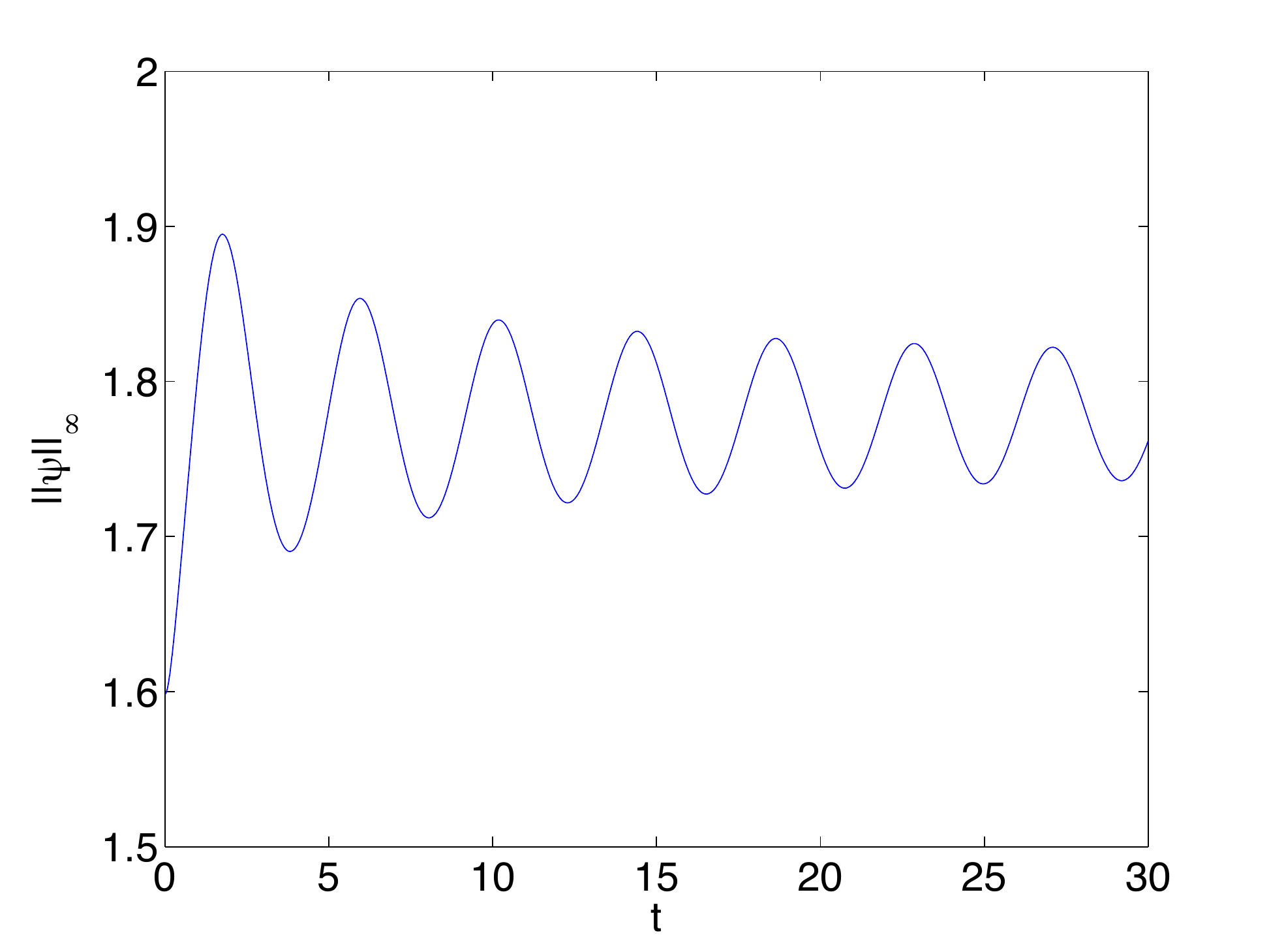}
 \caption{Time dependence of the $L^{\infty}$ norm of the solution of the focusing fNLS equation 
 (\ref{fNLS1}) with $s=0.9$ and initial data $\psi_{0}(x)=\alpha Q(x)$ for $\alpha=0.9$ (left) and $\alpha=1.1$ (right).}
 \label{sol_{09}_{09}max}
\end{figure}

\begin{figure}[htb!]
   \includegraphics[width=0.49\textwidth]{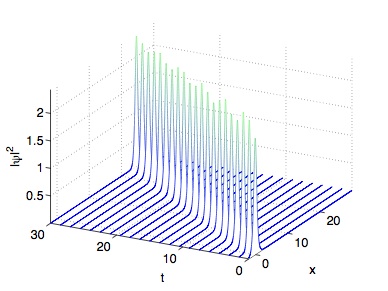}
   \includegraphics[width=0.49\textwidth]{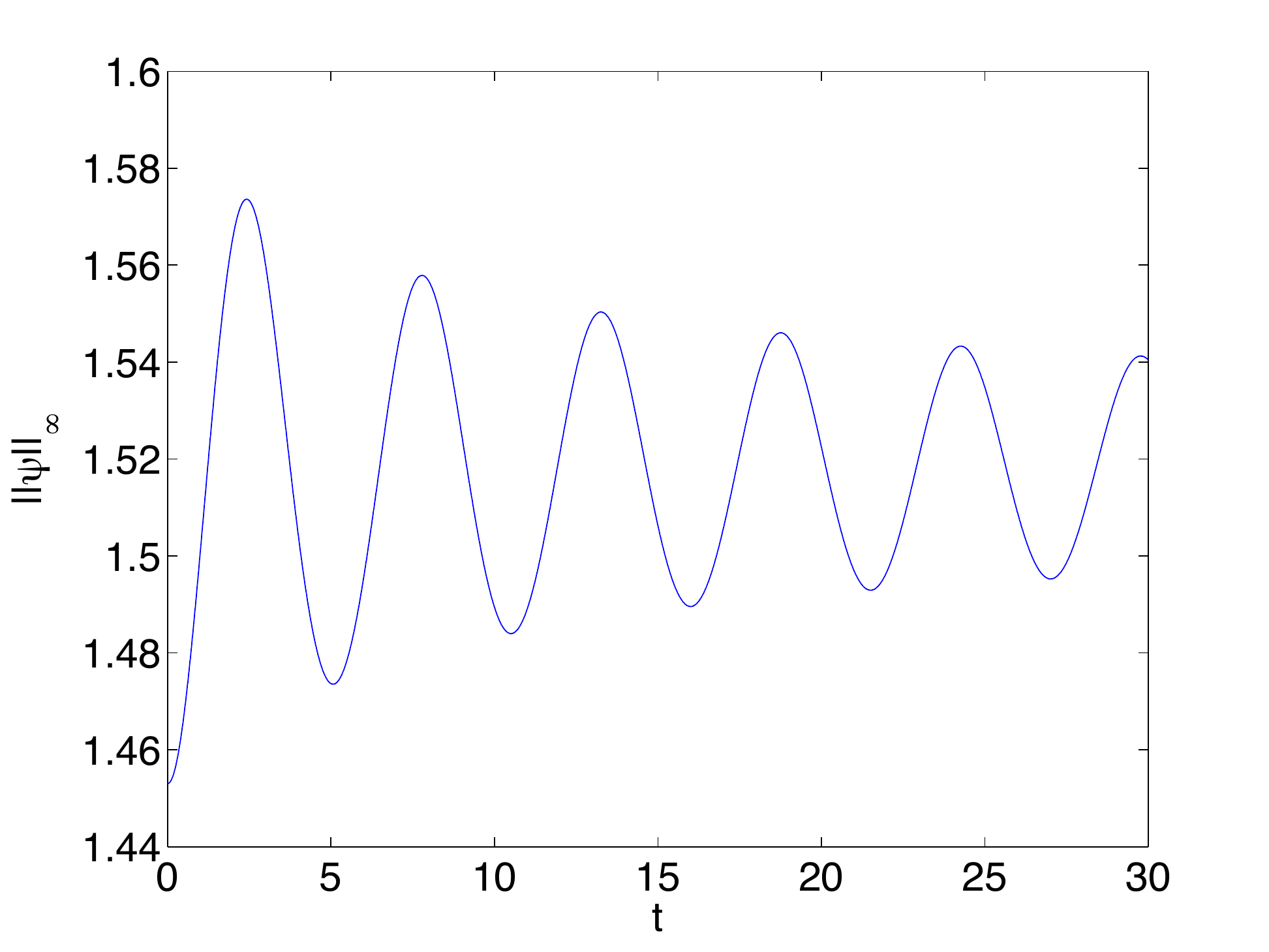}
 \caption{Time dependence of the modulus squared of the solution of the focusing fNLS equation 
 (\ref{fNLS1}) with $s=0.9$ and 
 initial data $\psi_{0}(x)=e^{ix} Q(x)$ (left). The 
behavior in time of the corresponding $L^{\infty}$ norm is given on the right.}
 \label{nlsffracsol_09_eix}
\end{figure}

For our final numerical test within this section, we first recall that the classical NLS equation ($s=1$) is {\it Galilei invariant}. This means, that if $\psi(t,x)$ is a solution, 
then so is 
\begin{equation}\label{eq:galilei}
\tilde{\psi}(t,x)=\psi(t,x-ct)e^{ic\cdot x +i|c|^{2}t/2},
\end{equation}
with $c\in \R^d$ some finite speed. In particular, if initially we choose $\psi(0,x)=Q(x)$, then we obtain the so-called {\it solitary wave} solution for NLS. 
For $s\not =1$, the Galilei symmetry of the model is broken, and hence we can not expect an exact formula of the same type as in \eqref{eq:galilei}.
Thus it is not obvious how an initial data of the form $\psi_0(x) = Q(x) e^{i x }$ (we set $c=1$ for simplicity) will evolve. However, it can be seen in Fig.~\ref{nlsffracsol_09_eix} that 
the initial hump still propagates essentially with constant velocity $\tilde c \approx 1$, similar to a solitary wave. 
The corresponding amplitude $|\psi(t,x)|$ oscillates around an asymptotically constant $L^{\infty}$ norm,  similar to the situation depicted in Fig.~\ref{sol_{09}_{09}max}. The latter is 
even more visible from the $L^{\infty}$ norm of the solution shown 
also in Fig.~\ref{nlsffracsol_09_eix}.
In other words, we find that initial data of the form $Q(x) e^{ix}$ 
give rise to a solution which can be seen as an {\it approximate solitary wave}, 
the amplitude of which converges as $ t\to +\infty$ to some, yet unknown, asymptotic profile.

\subsection{Perturbed ground states in the mass critical regime}

The mechanism described above, i.e., that a perturbed ground state 
asymptotically becomes a stationary state with the same mass as the 
initial datum, is not 
possible for the mass critical case $s=\frac{1}{2}$, since the $L^{2}$ norm 
and the equation are both invariant under the rescaling (\ref{scaling}). Thus it can be 
expected that the ground state is unstable in this case which is 
exactly what we observe for initial data of the form \eqref{iniQ}: First, for $\alpha=0.9$, i.e., initial data with mass smaller than the ground state, 
Fig.~\ref{sol_{05}_{09}} shows that 
the solution simply decays to zero with monotonically decreasing  
$L^{\infty}$ norm. Thereby the initial hump splits 
into two smaller humps which both move outwards.
\begin{figure}[htb!]
   \includegraphics[width=0.7\textwidth]{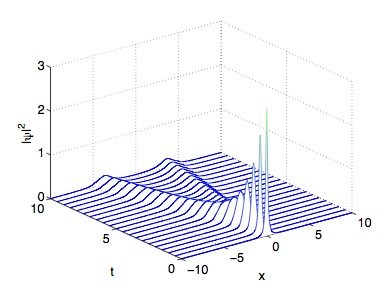}
 \caption{Modulus squared of the solution of the mass critical focusing fNLS equation 
 (\ref{fNLS1}) with $s=0.5$ and initial data $\psi_{0}(x)=0.9Q(x)$.}
 \label{sol_{05}_{09}}
\end{figure}

However, for an $\alpha>1$, i.e., a perturbation with mass larger than the 
ground state, the solution $\psi(t,x)$ appears to exhibit finite time blow-up, as 
can be seen in Fig.~\ref{sol_{05}_{11}}. \begin{figure}[htb!]
   \includegraphics[width=0.7\textwidth]{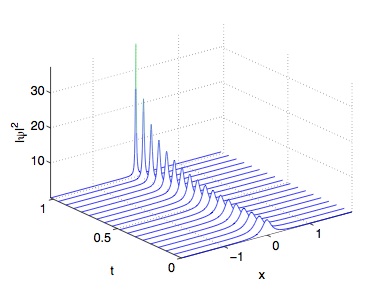}
 \caption{Modulus squared of the solution of the mass critical focusing fNLS equation 
 (\ref{fNLS1}) with $s=0.5$ and initial 
 data $\psi_{0}(x)=1.1 Q(x)$.}
 \label{sol_{05}_{11}}
\end{figure}
The blow-up is also 
indicated by the behavior of the $L^\infty$ norm and the $\dot H^1$ norm of the solution, see Fig.~\ref{sol_{05}_{11}norm}.
\begin{figure}[htb!]
   \includegraphics[width=0.49\textwidth]{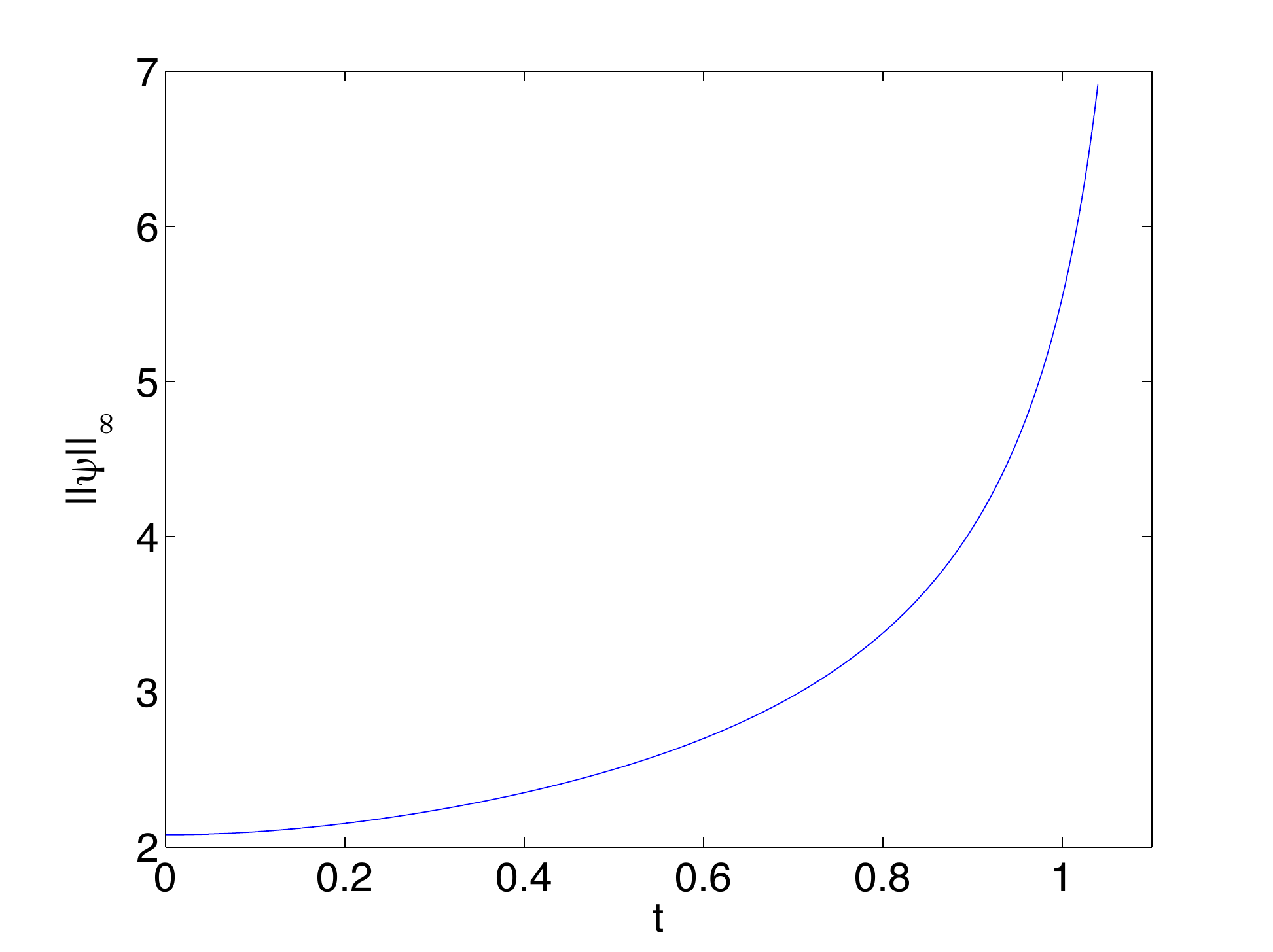}
   \includegraphics[width=0.49\textwidth]{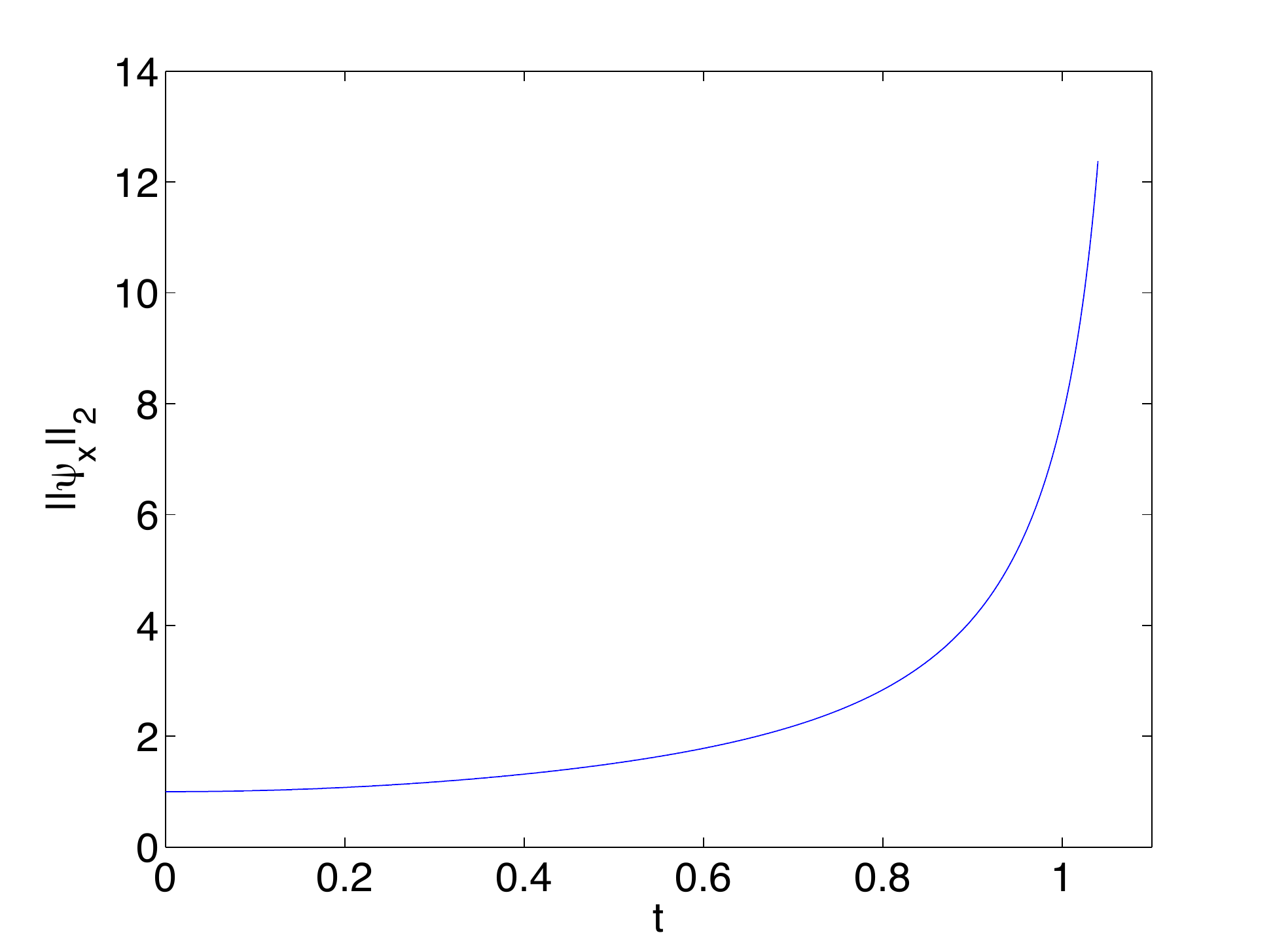}
 \caption{Time dependence of the $L^\infty$ norm (left) and of the 
 $\dot H^1$ norm normalized to 1 at $t=0$ (right) 
 of the solution of the focusing mass critical fNLS equation 
 (\ref{fNLS1}) with $s=0.5$ and initial 
 data $\psi_{0}(x)=1.1 Q(x)$.}
 \label{sol_{05}_{11}norm}
\end{figure}
Here, the Fourier coefficients are 
fitted to the asymptotic formula (\ref{fourasymp}). 
As explained above, the code is stopped once $\delta<m$, i.e., once the 
singularity is closer to the real axis than the smallest distance 
resolved by the Fourier method. 
Note that we run out of resolution in Fourier space before coming 
sufficiently close to the presumed blow-up. This is mainly due to the 
large computational domain $100[-\pi,\pi]$ which was needed because 
of the slow decrease of the ground state solution towards infinity. 
Around $t=1.0$ 
the resolution in Fourier space becomes insufficient, and the 
solution becomes incorrect as indicated by deterioration of the 
Fourier coefficients at a larger time. Note that the code would continue to 
run through and that the relative energy in this case would be 
still conserved to better 
than $10^{-9}$ at the final time. This shows that this quantity can 
only be used as an indicator if sufficient resolution in Fourier 
space is provided.


\section{Numerical studies of finite-time blow-up}\label{sec:blowup}

In this section, we will study the appearance of finite-time blow-up for 
solutions of the focusing cubic fNLS in $d=1$ 
with rapidly decreasing initial data $\psi_0\in \mathcal S(\R)$. The corresponding 
solution does not suffer from the same problems as the slowly decaying ground states, and hence gives 
a more reliable picture of the blow-up phenomena. 
Our choice of initial data is 
\begin{equation}\label{inisoliton}
\psi_0(x) = \frac{\beta }{\cosh( x)}\equiv \beta \, \text{sech}(x), \quad \beta \in \R,
\end{equation}
which are motivated by the soliton for the cubic NLS at $t=0$ (In particular these type of initial data have exponential decay as $|x|\to \infty$ which is preferable for our numerical studies).

\subsection{Numerical reproduction of known results for NLS} 

Before we investigate the blow-up for fNLS, we will test our numerical methods 
via a study of the focusing quintic NLS $p=2$ and septic NLS $p=3$. 
For the blow-up computations in this subsection, we always use 
$N=2^{17}$ Fourier modes for $x\in10[-\pi,\pi]$ and $N_{t}=50000$ 
time steps.

\begin{figure}[htb!]
   \includegraphics[width=0.49\textwidth]{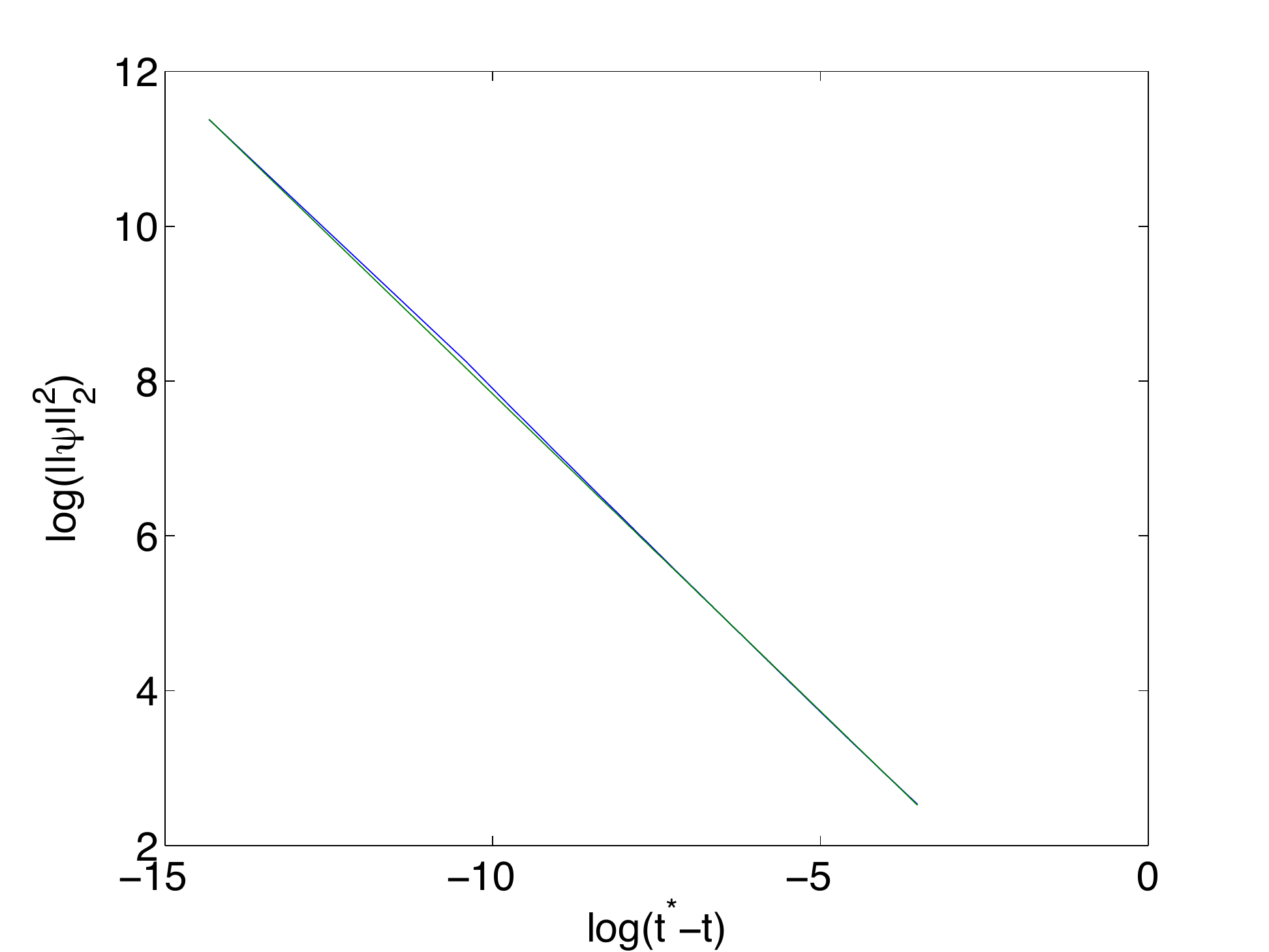}
   \includegraphics[width=0.49\textwidth]{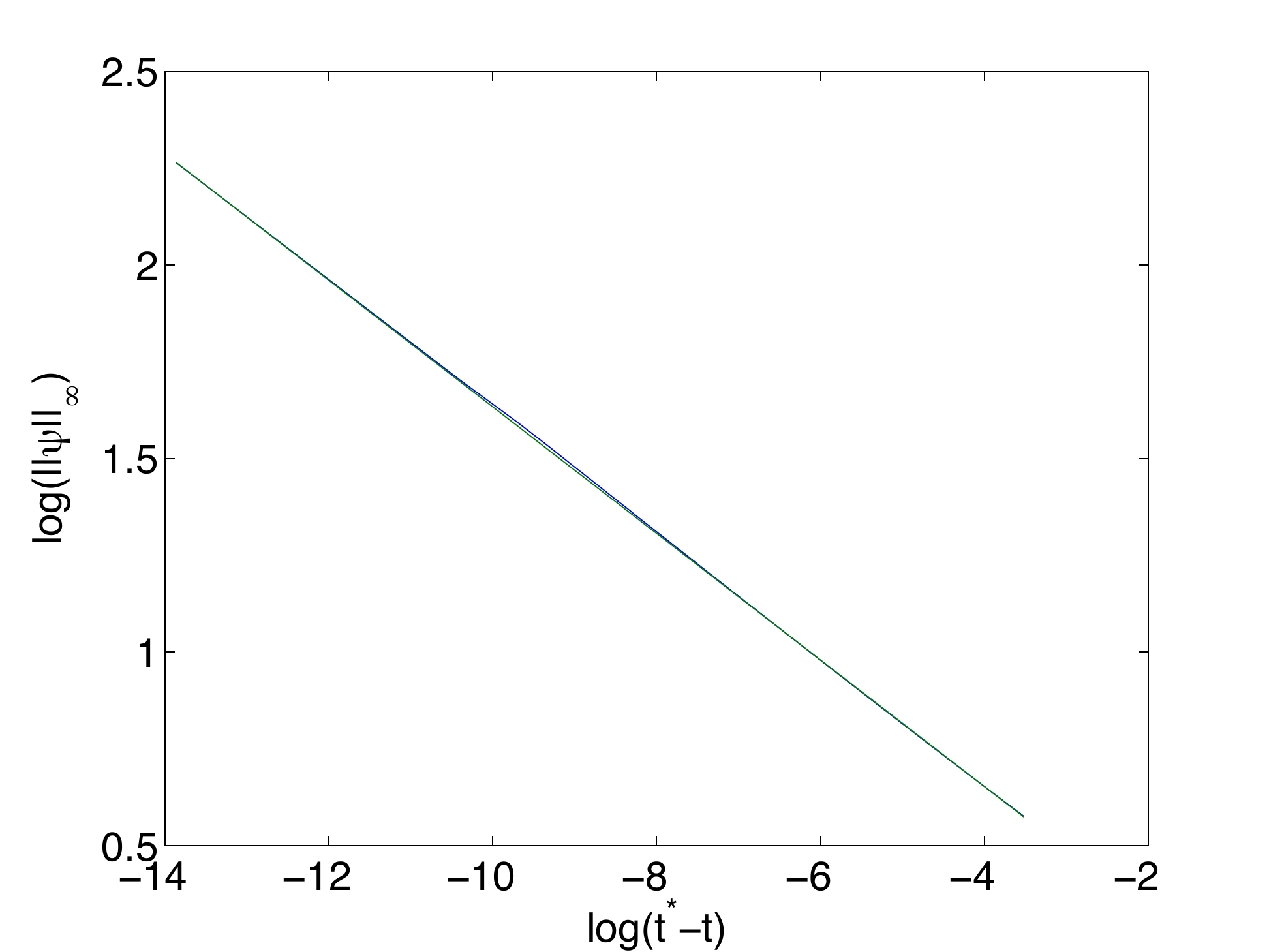}
 \caption{Fitting the logarithms of the $\dot H^1$ norm  
 (left) and of the 
 $L^{\infty}$ norm (right) of the solution of the septic 
 NLS equation ($s=1$) with initial data \eqref{inisoliton} close to the blow-up.  
The fitted line $\kappa_{1}\ln 
 (t^{*}-t)+\kappa_{2}$ (see the description) is given in green.}
 \label{nlsseptfig}
\end{figure}

We first consider the initial data \eqref{inisoliton} with $\beta = 1$ for the focusing 
septic NLS equation, i.e. \eqref{fNLS1} with $s=1$, $p=3$ and $\gamma = -1$. This 
equation is mass supercritical (and energy subcritical). We find that the 
code breaks at $t\approx1.4789$ due to an overflow 
error. The latter occurs in the computation of the nonlinearity 
$|\psi|^{2p}\psi$. At the last recorded time, the value of $\delta$ obtained by 
fitting the Fourier coefficients to the asymptotic formula 
(\ref{fourasymp}) is $\delta\approx2.4\times10^{-3}$ and thus more than an 
order of magnitude larger than the minimal resolved distance 
$m=4.794\times10^{-4}$ in (\ref{mres}). In order to obtain the actual blow-up time, we use 
the optimization algorithm \cite{fminsearch}, which is  
accessible via Matlab as the command \emph{fminsearch}. For $t\approx t^{*}$, we then fit for 
the $L^\infty$ and the $\dot H^1$ norm (always normalized to 1 at 
$t=0$ in this section) of $\psi(t,\cdot)$ to the expected asymptotic behavior 
(\ref{genscal}). The $L^{\infty}$ norm thereby catches the 
local behavior of the solution close to the blow-up point, whereas the 
homogenous Sobolev norm $\dot H^1$ 
takes into account the solution on the whole 
computational domain. Thus the consistency of the fitting results 
provides a test of the quality of the numerics. The results of the 
fitting can be seen in Fig.~\ref{nlsseptfig}. 
Fitting $\ln\|\partial_x \psi(t,\cdot) \|_{2}^{2}$ to $\kappa_{1}\ln (t^{*}-t)+\kappa_{2}$, we 
find $t^{*}=1.4789$, $\kappa_{1}=-0.8197\approx-5/6$ and $\kappa_{2}=-0.3644$. Similarly, 
we get for $\|\psi(t,\cdot) \|_{\infty}$ the values $t^{*}=1.4789$, 
$\kappa_{1}=-0.1634\approx -1/6$ and $\kappa_{2}=-0.0013$. Note the excellent 
agreement of the blow-up times which shows both the consistency of the 
fitting results and that the computation came very close to the 
blow-up. Note also the agreement with the predicted values $5/6$ 
respectively $1/6$ for the values of the $\kappa_{1}$. These values 
are unchanged within numerical precision if only the last 100 
computed time steps are used for the fitting.

The same initial data for the mass critical quintic NLS equation in $d=1$ 
lead to a breaking of the code at $t\approx 4.971$, again due to an overflow 
error. At the last recorded time, the value of $\delta$ obtained by 
fitting the Fourier coefficients to the asymptotic formula 
(\ref{fourasymp}) is $\delta\approx4.8\times10^{-3}$ and thus roughly an 
order of magnitude larger than the minimal resolved distance 
$m=4.794\times10^{-4}$. Fitting 
$\|\partial_x \psi(t,\cdot)\|_{2}^{2}$ as in the supercritical case 
to $\kappa_{1}\ln (t^{*}-t)+\kappa_{2}$, we 
get $t^{*}=4.9711$, $\kappa_{1}=-1.0077$ and $\kappa_{2}=1.3568$.
Similarly, we obtain for $\|\psi(t,\cdot)\|_{\infty}$ the values $t^{*}=4.9712$, 
$\kappa_{1}=-0.2533$ and $\kappa_{2}=0.3426$. The 
agreement of the blow-up times  shows again the consistency of the 
fitting results, and the agreement with the predicted values $1$ 
respectively $1/4$ for the values of the $\kappa_{1}$, if the scaling 
(\ref{genscal}) is assumed.

An important question is, whether the logarithmic corrections in (\ref{L2scal}) can also be seen within 
this approach. This is unlikely, since we do not use an adaptive 
rescaling here for the reasons explained before (that the periodic boundary conditions lead to numerical instabilities). 
To test what can be seen with the present code, we do the same fitting as 
above for the last 100 computed time steps since the logarithmic 
corrections will be mainly noticeable for $t\approx t^{*}$. In this case 
we get with numerical precision the same values for $\kappa_{1}$ and 
$\kappa_{2}$. We denote the $L^{2}$ norm of the difference between 
the logarithm of the fitted norm and $\kappa_{1}\ln 
(t^{*}-t)+\kappa_{2}$ as the \emph{fitting error} $\Delta_{2}$. We 
find $\Delta_{2}=1.88\times10^{-2}$ for the $L^2$ norm of $\psi_{x}$ and 
$\Delta_{2}=4.3\times10^{-3}$ for the $L^\infty$ norm of $\psi$. If we 
fit the same norms to $\tilde{\kappa}_{1}(\ln 
(t^{*}-t)-\ln\ln|\ln(t^{*}-t)|)+\tilde{\kappa}_{2}$, we get for the 
analogously defined fitting error $\tilde{\Delta}_{2}$ the values 
$3.72\times10^{-2}$ respectively $0.014$, i.e., higher values. Repeating 
the previous analysis for the last 10 computed  points, the fitting 
errors become $\Delta_{2}=8.7\times10^{-3}$ respectively $6.65\times10^{-2}$, 
and $\tilde{\Delta}_{2}=7.7\times10^{-4}$ respectively $9.9\times10^{-3}$, 
i.e., a better agreement for the logarithm corrections. Thus 
if there is an indication of the logarithmic corrections, they can 
only be expected very close to the time where the code is stopped for 
a lack of resolution. 

\subsection{Finite time blow-up for fNLS}

Having checked to which extent we are able to reproduce blow-up 
results for the usual NLS in $d=1$, we turn now to the case of \eqref{fNLS1} with $s\leq\frac{1}{2}$ 
and initial data given by \eqref{inisoliton}. By fitting the Fourier coefficients to the asymptotic formula (\ref{fourasymp}),
we find that a singularity is approaching the real axis in the complex plane 
for finite time, which indicates a blow-up. As discussed above, we 
stop the code once the value $\delta <m$ with $m$ given by (\ref{mres}). In contrast to the 
NLS examples with $p>1$, no overflow error is observed in the present 
case, due to the smaller power of the (cubic) nonlinearity. 
The blow-up time is again determined via 
the fitting of certain norms of the solution to the expected formulae 
(\ref{genscal}) and (\ref{L2scal}):

In the mass critical case $s=\frac{1}{2}$, the code is stopped at 
$t=2.9413$. Fitting, as before, the square of the $\dot H^1$ norm of $\psi$ for the 
last 1000 recorded time steps 
to $\kappa_{1}\ln (t^{*}-t)+\kappa_{2}$, we 
get $t^{*}=2.9940$, $\kappa_{1}=-2.0735$ and $\kappa_{2}=1.9783$. 
Similarly, we obtain for the $L^\infty$ norm of $\psi$ the values $t^{*}=2.994$, 
$\kappa_{1}=-0.5196$ and $\kappa_{2}=0.4003$. The 
agreement of the blow-up times  provides again a check of the consistency of the 
fitting results. The obtained values for $\kappa_{1}$ also agree well 
with the predicted values $-2$ 
respectively $-\frac{1}{2}$, if the scaling 
(\ref{genscal}) is assumed. To see whether there is an indication of 
logarithmic corrections to this formula as in (\ref{L2scal}), we 
repeat this fitting for the last 10 recorded time steps, just like in the 
$L^{2}$ critical case for $s=1$ above.  We 
find a fitting error $\Delta_{2}=4.5\times 10^{-3}$ for the $\dot H^1$ norm and 
$\Delta_{2}= 2.1\times10^{-3}$ for the $L^{\infty}$ norm of $\psi$. If we 
fit the same norms to $\tilde{\kappa}_{1}(\ln 
(t^{*}-t)-\ln\ln|\ln(t^{*}-t)|)+\tilde{\kappa}_{2}$, we get for the 
analogously defined fitting error $\tilde{\Delta}_{2}$ the values 
$4.8\times10^{-3}$ and $3.3\times10^{-5}$, respectively. In other words, we find essentially the same 
value for the $\dot H^1$ norm, but a much better agreement of the 
logarithmic correction for the $L^\infty$ norm of $\psi$ close to the blow-up. It is possible, however, 
that we did not get close enough to the blow-up in order for this effect 
to be also seen within the $\dot H^1$ norm, 
but it appears that the logarithmic correction is indeed visible locally near the blow-up. 

\begin{figure}[htb!]
   \includegraphics[width=0.49\textwidth]{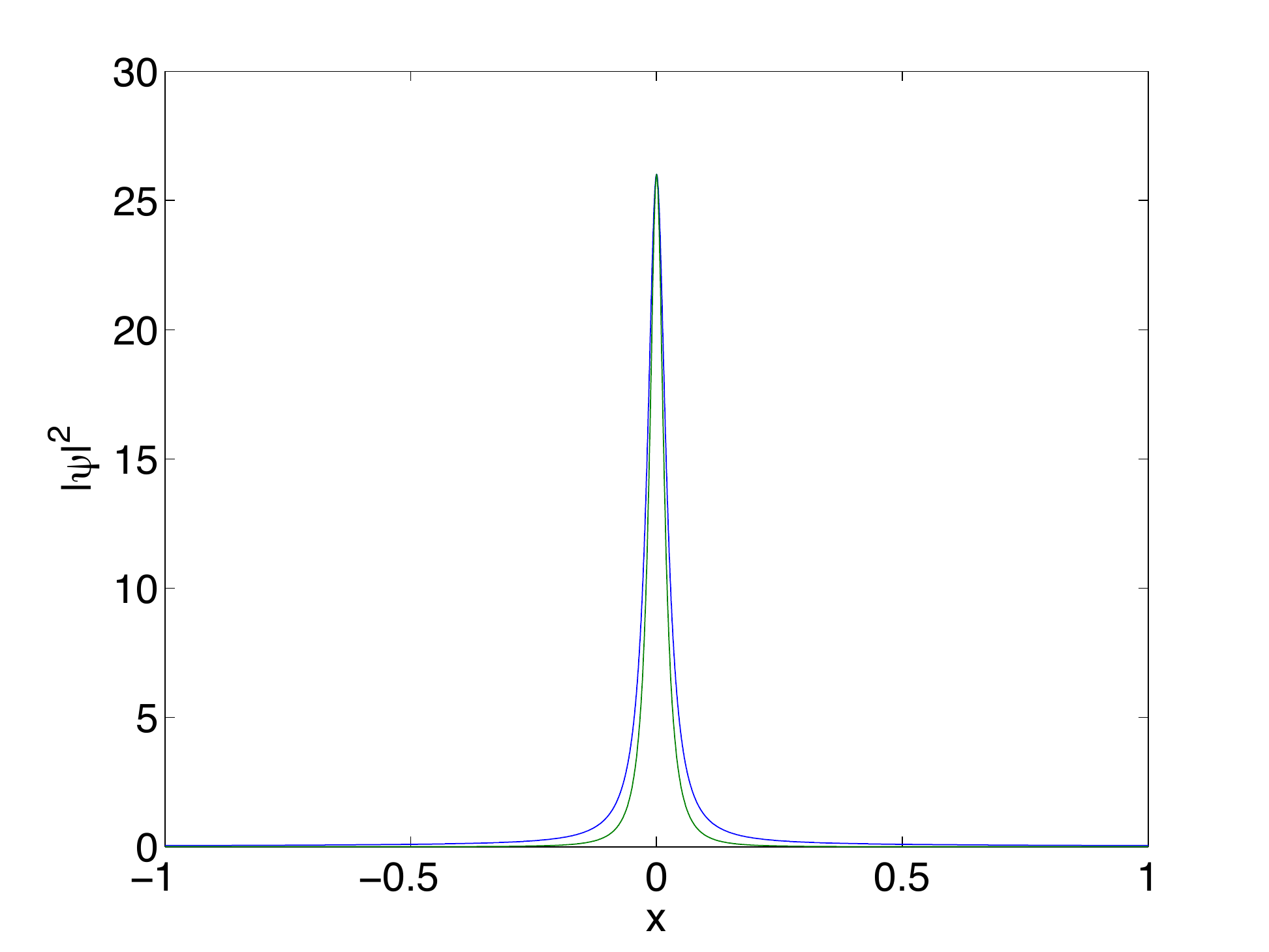}
   \includegraphics[width=0.49\textwidth]{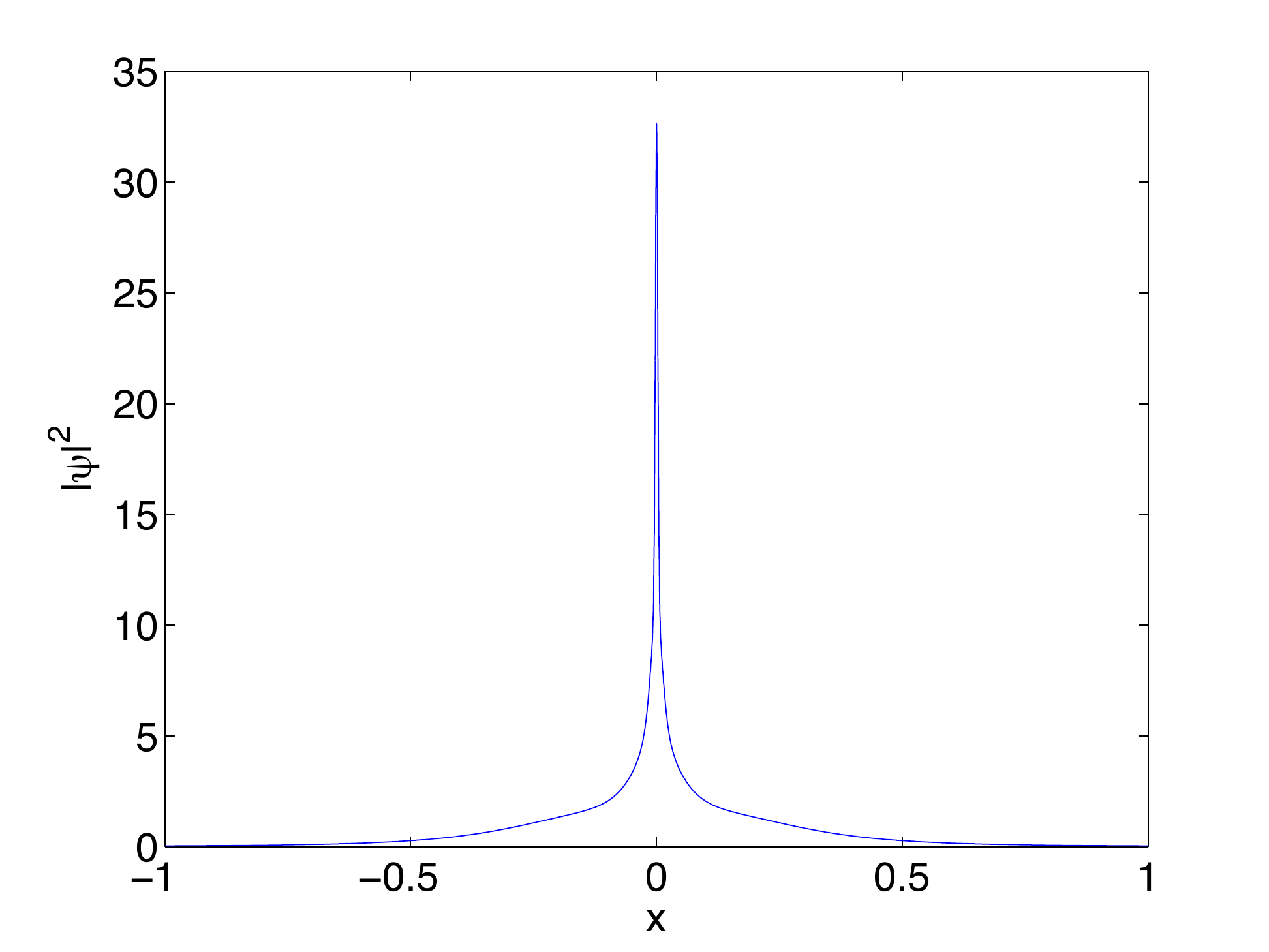}
 \caption{Blow-up profiles of solutions of the fNLS equation for the 
 initial data $\psi_{0}(x)=\mbox{sech}x$ at the respective
 last recorded time: on the left for the mass critical case $s=0.5$ 
 the modulus of the solution  in blue and the fitted soliton according 
 to the scaling (\ref{resc}); on the right for the energy 
 supercritical case $s=0.2$.}
 \label{blowupprofile}
\end{figure}

The blow-up profile of the solution at the last recorded time is 
shown in Fig.~\ref{blowupprofile} on the left. In the same figure we 
show the soliton  rescaled according to (\ref{resc}). The scaling 
factor $L$ is simply fixed in a way that the maxima of the solutions coincide. 
It can be seen that the agreement is qualitatively good and 
quantitatively convincing close to the maximum. Obviously the 
asymptotic description is less satisfactory for a larger distance to 
the maximum due to the slow algebraic decrease of the soliton for 
$|x|\to\infty$. It is to be expected that the asymptotic description 
would improve if times closer to blow-up could be reached. 

In the mass supercritical case $s=0.4$, we observe for the same 
initial data that the code is stopped at a larger time $t=3.1347$ than in the case $s=\frac{1}{2}$.
Once again we fit the $\dot H^1$ norm of $\psi$ for the 
last 1000  recorded time steps 
to $\kappa_{1}\ln (t^{*}-t)+\kappa_{2}$ and  
obtain $t^{*}=3.1396$, $\kappa_{1}=-2.3521$ and $\kappa_{2}=2.5182$; 
similarly we obtain for the $L^\infty$ norm of $\psi$ the values $t^{*}=3.1396$, 
$\kappa_{1}=-0.5192$ and $\kappa_{2}=0.4441$. There is again a good 
agreement of the fitted blow-up times  and of the values of the $\kappa_{1}$ 
with the predicted values $-2.25$ and $-\frac{1}{2}$, respectively, if the scaling 
(\ref{genscal}) is assumed.

\begin{figure}[htb!]
   \includegraphics[width=0.7\textwidth]{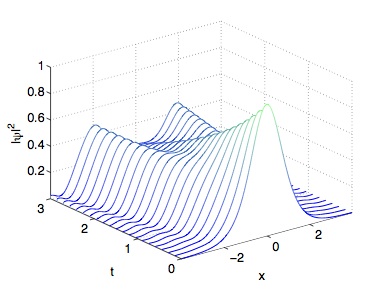}
 \caption{Modulus squared of the solution to the mass supercritical focusing fNLS equation 
 (\ref{fNLS1}) with $s=0.4$ and initial data 
 $\psi_{0}(x)=e^{ix^{2}}\mbox{sech}(x)$.}
 \label{nlsfsechs04_eix2}
\end{figure}

It is known (see, e.g., \cite{Caz}) that multiplication of the 
initial data with a rapidly oscillating factor of the form $e^{ib|x|^{2}}$ with $b>0$ 
introduces a defocusing effect in the standard focusing NLS 
equation. Indeed, one can prove that for $b>0$ sufficiently large, the 
solution of 
NLS exists for all $t\ge0$, regardless of the sign of $\gamma$ in 
(\ref{fNLS1}).
Again these analytical considerations do not directly 
carry over in the presence of fractional derivatives. But it can be seen in Fig.~\ref{nlsfsechs04_eix2} that 
the fNLS solution for $s=0.4$ and $\psi(0,x)=e^{ix^{2}}\mbox{sech}(x)$ 
not only does {\it not} show blow-up as above, but displays the behavior of 
solutions to the defocusing fNLS equation to be studied in the  
sections below.


\section{The energy critical and supercritical regime}\label{sec:energy}

Recall that there is no energy supercritical regime for the usual NLS in $d=1$. 
However, we can reach this regime in the fractional NLS \eqref{fNLS1}  as soon as $s<\frac{1}{4}$. 
We shall study both, the focusing and the defocusing situations in more 
detail. For cases where no blow-up is observed, we also trace the 
$\dot H^\sigma(\R)$ norm invariant under the rescaling (\ref{resc}). To this end, we
consider here $s=0.2$ for which the critical exponent is $\sigma = 0.3$, in view of (\ref{dotHexp}). 
The initial data will be the same as in Section \ref{sec:blowup} above, i.e.,
\[
\psi_0(x) =  \beta \, \text{sech}(x), \quad \beta \in \R.
\]

\subsection{Finite time blow-up for energy supercritical fNLS}
\begin{figure}[htb!]
   \includegraphics[width=0.49\textwidth]{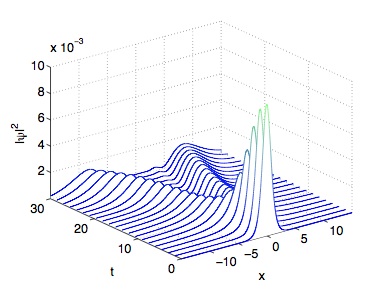}
   \includegraphics[width=0.49\textwidth]{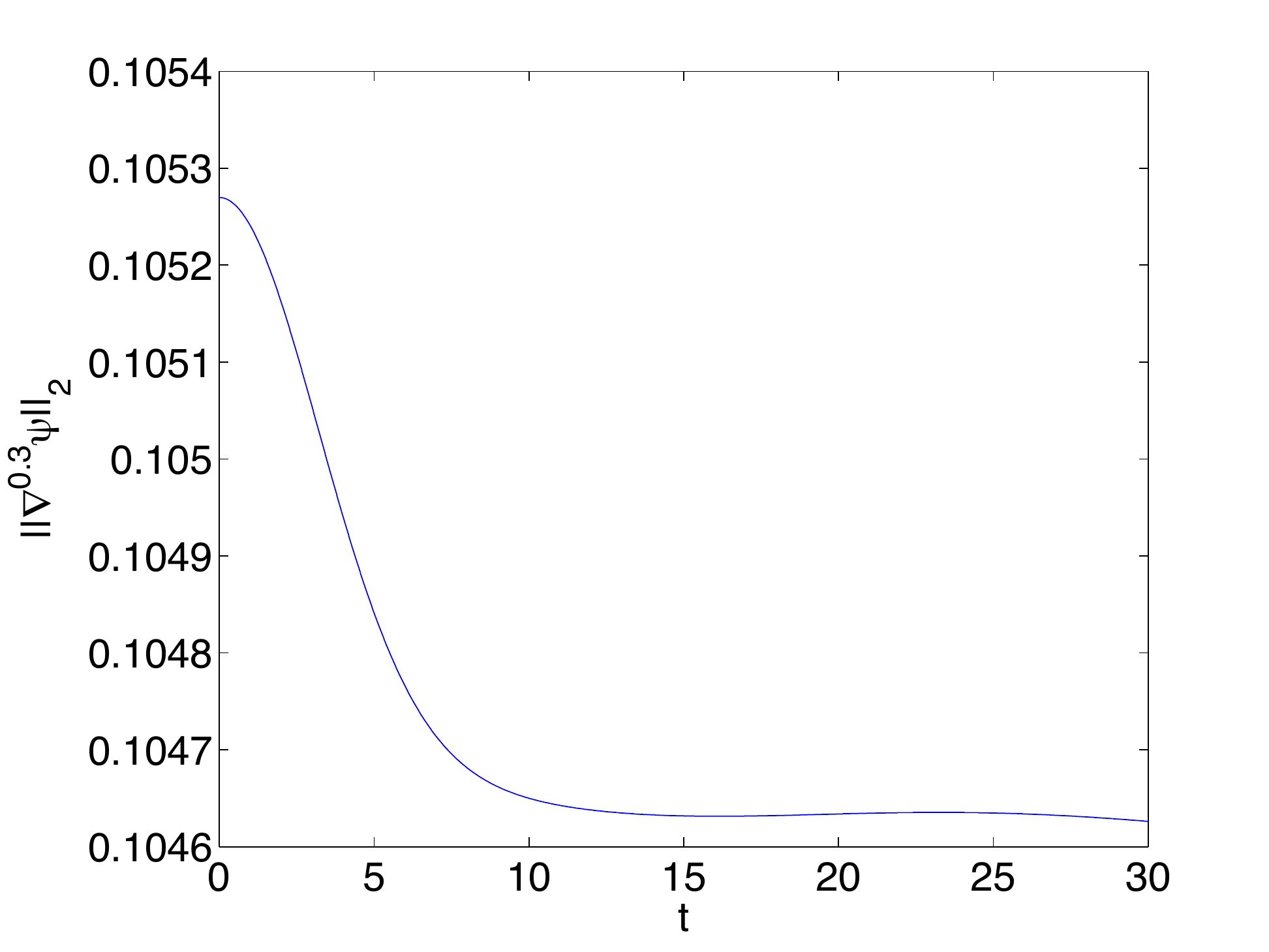}
 \caption{Solution of the energy supercritical focusing fNLS equation 
 (\ref{fNLS1}) with $s=0.2$ and initial 
 data $\psi_{0}(x)=0.1 \mbox{sech}( x)$ on the left, and the associated invariant $\dot 
 H^{0.3}$ norm of the solution on the right.}
 \label{nlsfrac01sechs02}
\end{figure}

It is not clear what the precise conditions on the 
initial data are, which lead to finite-time blow-up in the energy supercritical regime. The 
numerical experiments of the previous subsection seem to indicate 
that initial data in the vicinity of the ground state with larger 
mass and smaller energy than the ground state produce such a blow-up. 
In fact, if we study initial data with small 
mass, say $\psi_{0}=0.1\mbox{sech}(x)$, we find that the initial hump simply decays to zero as $t\to +\infty$, as can be seen in 
Fig.~\ref{nlsfrac01sechs02}. The $L^{\infty}$ norm of the solution is 
monotonically decreasing and there is no indication of blow-up in this 
case. The scaling invariant $\dot H^{0.3}$ norm also appears to be bounded as 
can be seen in the same figure.

However, for the initial data $\psi_{0}=\mbox{sech} (x)$, the code is 
stopped at the time $t=6.0748$ since the distance between 
a singularity (as indicated by 
the Fourier coefficients via (\ref{fourasymp})) and the real 
axis is smaller than the numerical resolution. Fitting, as before, the norm $\dot H^1$ norm of $\psi$ for the 
last 1000 time steps 
to $\kappa_{1}\ln (t^{*}-t)+\kappa_{2}$, we 
get $t^{*}=6.2771$, $\kappa_{1}=-3.6949$ and $\kappa_{2}=8.5231$ with 
a fitting error $\Delta_{2}\approx 10^{-2}$.
Similarly we obtain for the $L^\infty$ norm of $\psi$ the values $t^{*}=6.2804$, 
$\kappa_{1}=-0.5779$ and $\kappa_{2}=1.2001$ with a fitting error 
$\Delta_{2}$ of 
the order of $10^{-3}$. These values for $\kappa_{1}$ agree 
with the predicted values $-3.5$ 
respectively $-\frac{1}{2}$, if the scaling 
(\ref{genscal}) is assumed. Note that the blow-up time $t^{*}$ is 
more than {\it twice} the $t^{*}$
found for the same initial data in the mass critical case, which seems quite surprising (as one would naively expect the blow-up time to be monotonically dependent on the choice of $s$). 
The agreement of the fitted blow-up times for the two norms is worse
than in the mass critical case. This is due to the fact that the code did not get as 
close to the blow-up time as for $s=\frac{1}{2}$. It appears that 
considerably higher resolution would be needed in this case as is 
indicated by the stronger divergence of the $\dot H^1$ norm. In Fig.~\ref{blowupprofile} we show on 
the right the blow-up profile of the solution at the last recorded 
time. Visibly this profile is different from the blow-up in the mass 
critical case which qualitatively corresponds to the soliton. Here 
the blow-up profile is much more compressed which also explain why we 
could not get as close to the blow-up as in the mass critical case.

\subsection{Long time behavior for defocusing energy supercritical fNLS}

\begin{figure}[htb!]
   \includegraphics[width=0.49\textwidth]{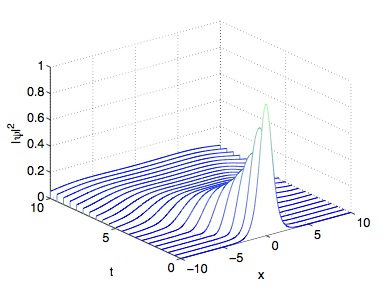}
   \includegraphics[width=0.49\textwidth]{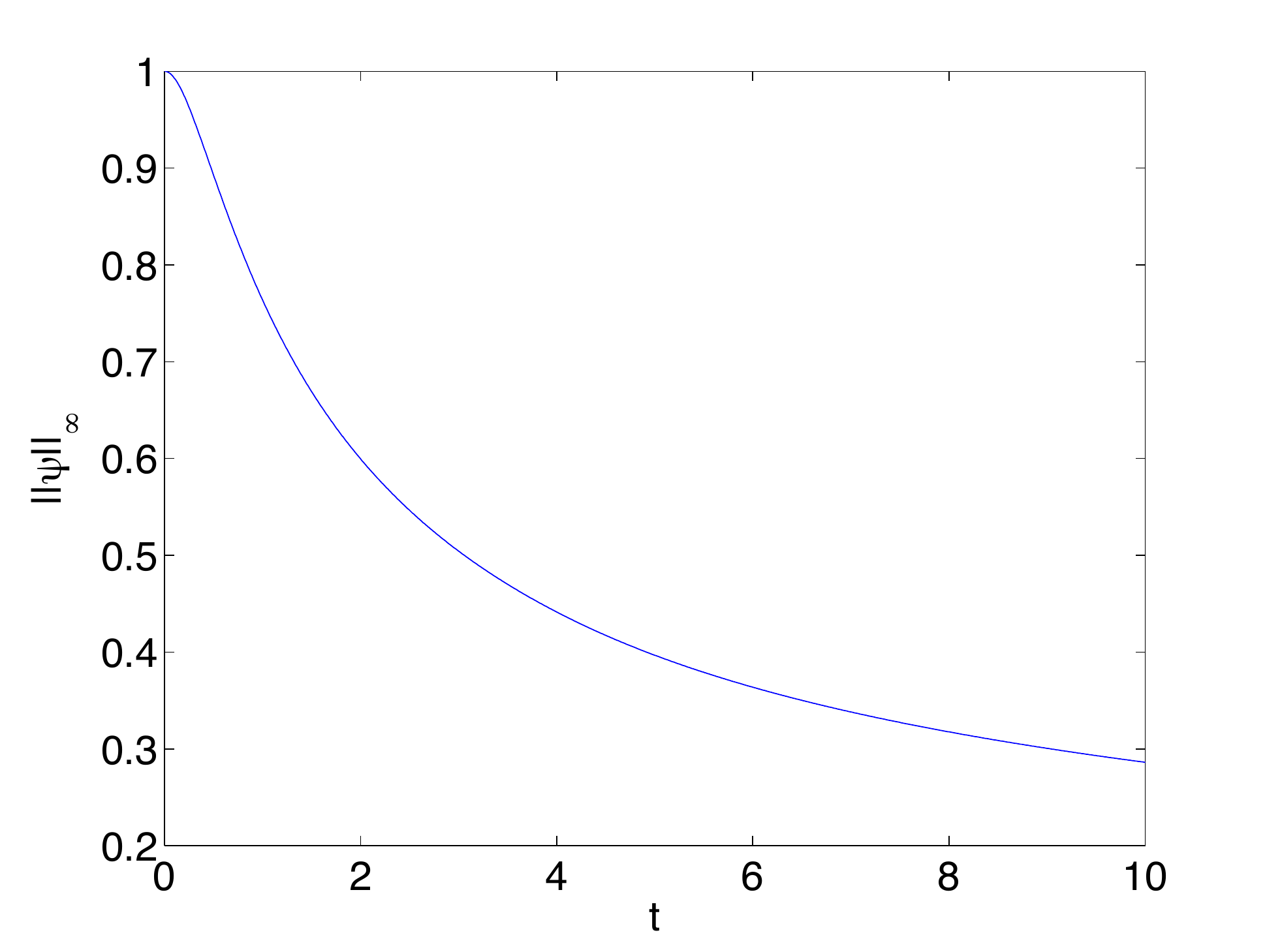}
 \caption{Left: Solution of the mass subcritical defocusing fNLS equation (\ref{fNLS1}) with $s=0.9$ and initial 
 data $\psi_0= \text{sech}(x)$. Right: The corresponding 
 $L^{\infty}$ norm of the solution.}
 \label{nlsfracdsech09}
\end{figure}

In Fig.~\ref{nlsfracdsech09} we show the solution of a defocusing fNLS 
equation (\ref{fNLS1}) with $s=0.9$ and initial data (\ref{inisoliton}) with $\beta=1$. It can be 
seen that the time evolution of the solution simply disperses the 
initial datum. 
This behavior is even more visible from the $L^\infty$ norm of the 
solution which is shown in the same figure. Obviously 
the norm is monotonically decreasing. 

\begin{figure}[htb!]
   \includegraphics[width=0.7\textwidth]{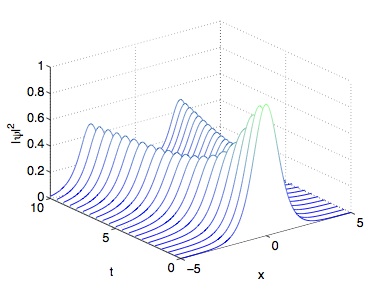}
 \caption{Solution of the defocusing energy supercritical fNLS equation (\ref{fNLS1}) with $s=0.2$ and initial 
 data $\psi_0= \text{sech}(x)$.}
 \label{nlsfracdsech02}
\end{figure}

For weaker dispersion, i.e., smaller $s$, the situation changes the 
shape, 
however. As can be seen in  
Fig.~\ref{nlsfracdsech02}, for $s=0.2<\frac{1}{4}$, i.e. the energy supercritical regime, the initial hump splits into two humps 
both of which travel to spatial infinity.  

\begin{figure}[htb!]
   \includegraphics[width=0.45\textwidth]{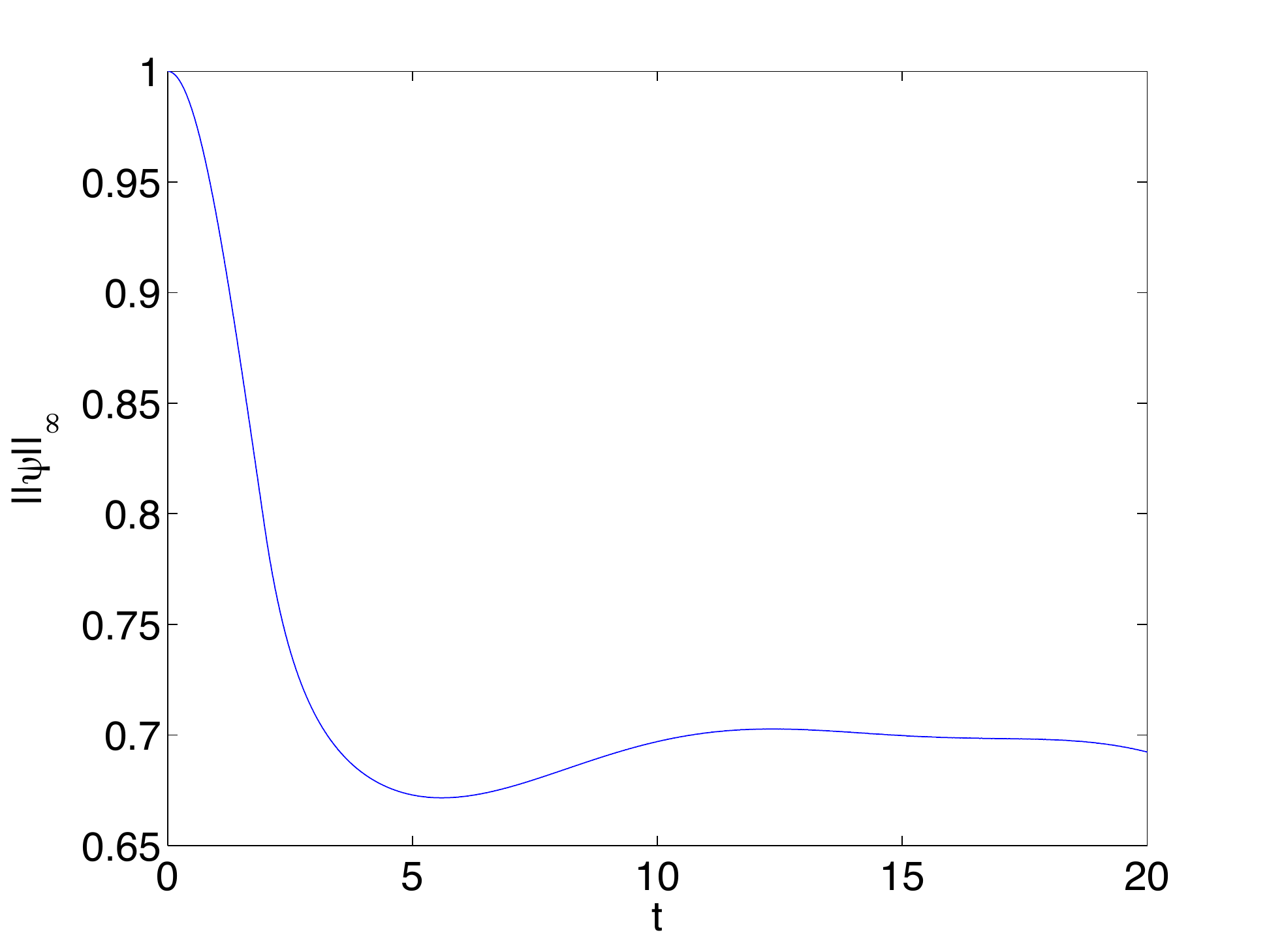}
   \includegraphics[width=0.45\textwidth]{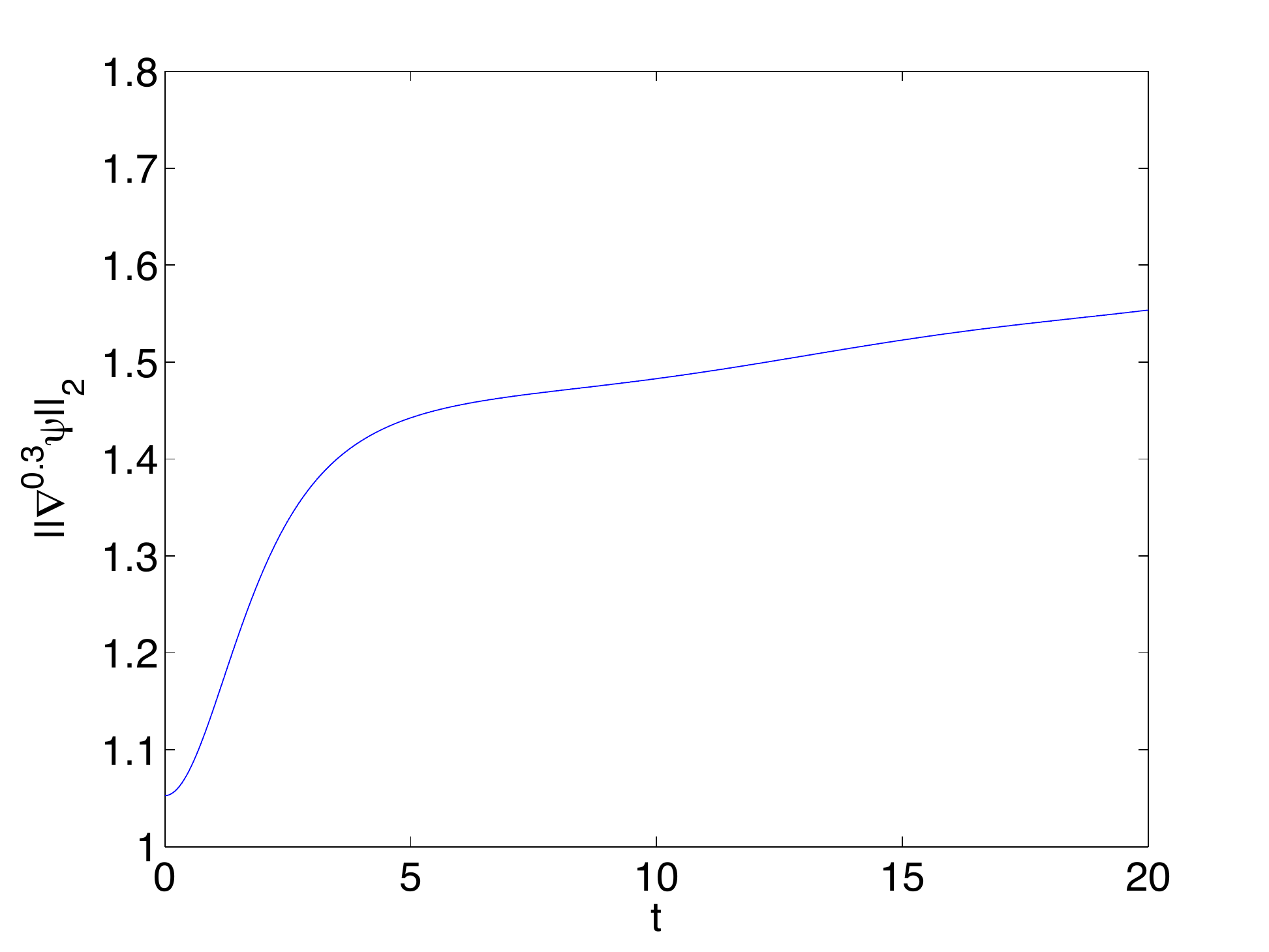}
 \caption{Time dependence of the $L^\infty$ norm (left) and of the 
 scaling invariant $\dot H^{0.3}$ norm (right) of the solution of the energy supercritical 
 defocusing fNLS equation 
 (\ref{fNLS1}) with $s=0.2$ and initial 
 data $\psi_0= \text{sech}(x)$.}
 \label{nlsfracdsech02max}
\end{figure}

The formation of the two humps can also be inferred from the $L^\infty$ norm of the 
solution which is shown in Fig.~\ref{nlsfracdsech02max}. 
The $L^\infty$ norm 
appears to become almost constant for large $t\gg 0$. An interesting 
quantity in this context is the scaling invariant $H^{0.3}$ norm
since its boundedness would allow to control the solution globally in 
time. It can be seen that there is no indication that this norm 
diverges (also not for larger times than shown here). However, in comparison to the earlier numerical study of \cite{CSS} for energy supercritical NLS (in $d=5$), we do not find 
that the critical $H^{0.3}$ norm approaches a constant for $t\gg 0$. This is probably due to the fact that in our case, the solution 
decays much slower as $x\to \infty$ which prevents our numerical 
method from computing for sufficiently long times. A possible way to overcome this issue is the scaling regime introduced in the 
next section.


\section{Numerical study of the semiclassical regime}\label{sec:semi}

\subsection{Semiclassical rescaling} 

A possible approach to study the long time behavior of 
solutions of the (dimensionless) fNLS (\ref{fNLS1}) is by considering slowly varying initial data of the form
\[
\psi_0(x;\epsilon)=  u(\epsilon x), 
\]
where $0<\epsilon \ll 1$ is a small {\it semiclassical parameter} and $u\in \mathcal S(\R^d)$ is some given initial profile. As $\epsilon \to 0$ the 
initial data approaches the constant value $u(0)$. Hence, in order to see nontrivial effects one has to wait until sufficiently long times of order $t\sim O(1/\epsilon)$, which 
consequently requires to rescale the spatial variable onto macroscopically large scales $x\sim O(1/\epsilon)$ too. 
In other words, we consider $x\mapsto \tilde x= x\epsilon$, $t\mapsto \tilde t= t\epsilon$ and set 
\[
\psi^\e(\tilde t,\tilde x) = \psi(\tilde t/\epsilon, \tilde x/ \epsilon),\] 
to obtain the following {\it semiclassically scaled} fNLS for the new unknown $\psi^\e$, where 
we discard the ``tildes" again for the sake of simplicity:
\begin{equation}
  i \epsilon\partial_t \psi^\epsilon = \frac{\epsilon^{2s}}{2}{(-\Delta)}^s \psi^\e + \gamma |\psi^\e|^{2} \psi^\e, \quad \psi^\e(0,x)=u(x),\quad x\in \R^d. 
    \label{fNLSe}
\end{equation}
For $\epsilon\ll 1$, the behavior of this equation describes solutions on macroscopically large space and time-scales, which justifies the name of semiclassical asymptotics. 
From the numerical point of view, the simulation of Schrödinger equations in the semiclassical regime is a formidable challenge, 
for which several different techniques have been developed in recent years, see the review \cite{JMS} for more details.


For the classical cubic NLS ($s=1$), an asymptotic theory for the solution as $\e \to 0$ is usually based on WKB type expansions. To this end, one assumes
\begin{equation}\label{eq:ansatz}
\psi(t,x) \approx a(t,x)e^{i S(t,x)/\epsilon},
\end{equation}
where $S(t,x)\in \R$ is a real-valued phase function an $a(t,x)\in \mathbb C$ a (in general) complex-valued amplitude. In the defocusing case, one can prove (see, e.g. \cite{Carles} and the references therein) 
that as $\epsilon \to 0$, this gives a valid approximation of the exact solution $\psi^\e$ provided $a, S$ are sufficiently smooth solutions of the following hydrodynamic system:
\begin{equation}\label{WKB}
\begin{cases}\partial_t S + \frac{1}{2} |\nabla S|^2 + \gamma |a|^2 = 0, \\
\partial_t a + \nabla S \cdot \nabla a + \frac{a}{2} \Delta S = 0, \end{cases}
\end{equation}
or, in terms of $\rho=|a|^2 $ and $v=\nabla S$:
\begin{equation}\label{qhd}
\begin{cases}\partial_t v + v\cdot \nabla v +  \gamma \nabla \rho =0, \\
\partial_t \rho + \text{div} (v\rho) = 0. \end{cases}
\end{equation}
Since this system in general exhibits shocks, the WKB approximation is only valid for short times $t< t_c$, where $t_c>0$ is the time where the first 
shock appears (also known as caustic-onset time, or time of the first gradient catastrophe). In the focusing case, the situation is even worse, as the obtained hydrodynamic system is 
found to be elliptic and thus not well-posed, see  \cite{KMM,TVZ} for partial results in the completely integrable case $d=p=1$.
\begin{remark}
For $s=1$ the system \eqref{WKB} can be formally obtained from the so-called Madelung system in the limit $\e \to 0$. The latter is obtained by inserting 
the right hand side of \eqref{eq:ansatz} into the NLS and separating real and imaginary parts, which gives
\begin{equation}\label{WKB}
\begin{cases}\partial_t S + \frac{1}{2} |\nabla S|^2 + \gamma |a|^2 = 0, \\
\partial_t a + \nabla S \cdot \nabla a + \frac{a}{2} \Delta S = \epsilon^2 \frac{\Delta a}{a}. \end{cases}
\end{equation}
This system is indeed equivalent to the NLS,  provided $a\not =0$.  This formulation has been used in, e.g., \cite{JLM} to study the semiclassical limit of defocusing NLS. 
In the case of fractional NLS, no such equivalent Madelung type equivalent system has yet been derived (due to the lack of an appropriate Leibnitz rule).
\end{remark}

\subsection{Semiclassical limit of the focusing fractional NLS equation}
\begin{figure}[htb!]
   \includegraphics[width=0.49\textwidth]{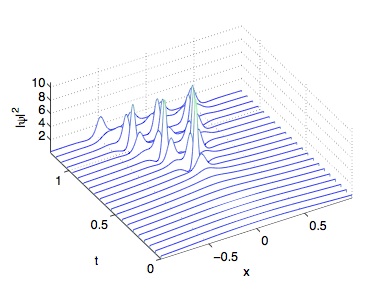}
   \includegraphics[width=0.49\textwidth]{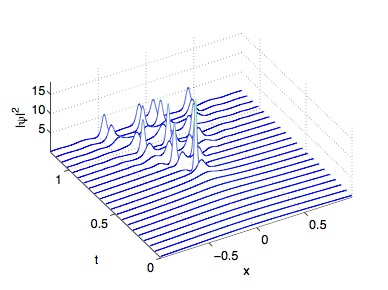}
  \caption{Solution of the semiclassical, focusing fNLS equation (\ref{fNLSe})
 for $\epsilon=0.1$ and the initial data 
 $\psi_{0}=\mbox{sech}(x)$; on the left for $s=1$, on the 
 right for $s=0.9$.}
 \label{nlsfseche01s09}
\end{figure}

In this subsection we study the semiclassically scaled, focusing (cubic) fNLS equation (\ref{fNLSe}) 
with $\gamma=1$ and $\e \ll 1$. We choose 
$\psi^{\epsilon}(x,0)=\mbox{sech}x$. For NLS ($s=1$) in $d=1$ and generic initial data, it is known that the semiclassical 
system \eqref{qhd} exhibits a gradient catastrophe at some finite time $0<t_c<+\infty$, yielding a square root type singularity 
in the gradient of the phase, see \cite{DGK}. 
For $\epsilon>0$, this singularity is regularized by highly oscillatory waves (the so-called dispersive shock phenomena). Indeed, one can see numerically that 
the square-modulus of the solution continues to grow for some time $t>t_{c}$, and eventually splits into several smaller humps leading to a zone of 
modulated oscillations as can be seen in Fig.~\ref{nlsfseche01s09}. 

\begin{remark}We refer to \cite{DGK}, \cite{DGK2} for more details and a conjecture concerning the asymptotic description of 
the solution of semiclassical NLS near $t\sim t_{c}$. See also \cite{BT} for a partial proof of this conjecture. \end{remark}

For smaller values of $s<1$ and $\epsilon<1$, the dispersion gets weaker 
and the focusing effect of the equation becomes stronger. This leads to a 
higher maximum and a more ``agitated'' oscillatory zone after the maximal 
peak. It was argued in \cite{etna}, that this maximal peak needs to be 
numerically well resolved. If this is not the case, the Fourier 
coefficients for the high wave numbers get polluted, which triggers 
the \emph{modulation instability} of the focusing NLS equation. The latter phenomenon 
cannot be controlled even with Fourier filtering methods.  Therefore we cannot 
reach much smaller values for $s$ and $\epsilon$ than used below. It would be 
necessary to go to higher than double precision to be able to address 
more extreme 
cases. In Fig.~\ref{nlsfseche01s09} we show the solution of the 
semiclassical focusing fNLS for $\epsilon=0.1$, $s=0.9$, and initial data 
 $\psi^\e(0,x)\equiv u(x) = \text{sech}(x)$.  The 
computation is carried out here with $N=2^{16}$ Fourier modes and 
$N_{t}=20000$ time steps. 

\begin{figure}[htb!]
   \includegraphics[width=0.49\textwidth]{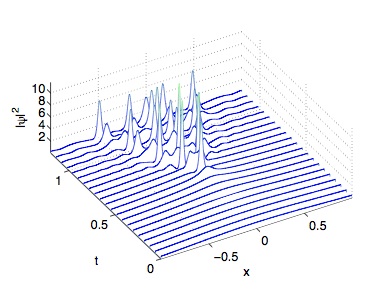}
   \includegraphics[width=0.49\textwidth]{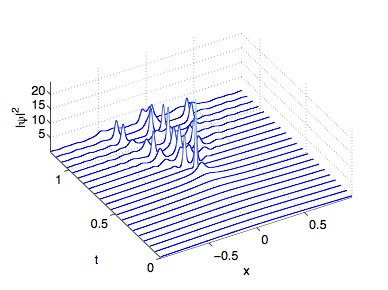}
  \caption{Solution of the semiclassical, focusing fNLS equation (\ref{fNLSe})
with initial data $\psi^\e_{0}=\mbox{sech}(x)$; on the left for $\epsilon=0.08$ and 
 $s=0.9$, on the right for $\epsilon=0.1$ and $s=0.08$.}
 \label{nlsfseche01s08}
\end{figure}

The same resolution can be used to study the same situation for the 
slightly smaller value of $\epsilon=0.08$. It can be seen in 
Fig.~\ref{nlsfseche01s08} that, as expected, there are more 
oscillations and a higher maximum in this case. To treat the solution for the same initial data with even smaller 
$s=0.8$, we had to use $N=2^{18}$ Fourier modes and $N_{t}=50000$ 
time steps in 
Fig.~\ref{nlsfseche01s08}. It is clearly visible that the maximum of 
the solution continues to grow as expected with smaller $s$, and that 
the oscillatory zone shows more humps than for the same value of 
$\epsilon$, but larger $s$. As discussed in the previous sections, a 
blow-up is to be expected for sufficiently small $\epsilon$ for 
$s\leq \frac{1}{2}$.

\subsection{Semiclassical limit of the defocusing fractional NLS equation}
\begin{figure}[htb!]
   \includegraphics[width=0.49\textwidth]{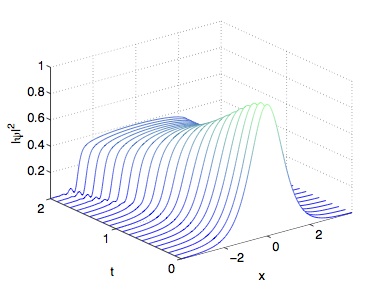}
   \includegraphics[width=0.49\textwidth]{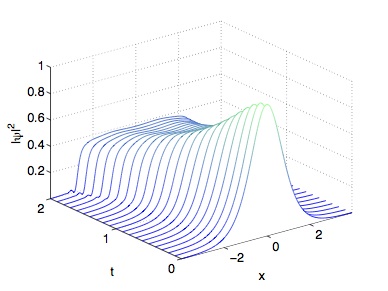}
  \caption{Solution of the semiclassical, defocusing fNLS equation (\ref{fNLSe})
with $\epsilon=0.1$, and initial data 
 $\psi_{0}=\mbox{sech}(x)$, on the left for $s=1$, on the right for 
 $s=0.9$.}
 \label{nlsfdseche01s09}
\end{figure}

In this subsection we study the the semiclassically regime for defocusing fNLS 
equations. In the NLS case ($s=1$), it is known that solutions corresponding to initial data  $ \psi^\e(0,x)\equiv u(x) = \text{sech}(x)$ exhibit 
a gradient catastrophe at two points, here, for symmetry reasons, at 
$\pm x_{c}$. As in the case of solutions of the Hopf equation, this 
is a cubic singularity at the onset of the formation of a shock. For 
small $\epsilon>0$, this singularity is regularized in the form 
of a zone of rapid modulated oscillations as in dispersive shocks of 
the KdV equation. The 
initial hump is defocused whilst the sides of the hump steepen. At a 
given point, oscillations form near these strong gradients. In 
Fig.~\ref{nlsfdseche01s09}, small oscilations  appear near the hump 
on the left (and hidden by the hump on the right) at the last shown 
times.

For smaller values of $s$ and $\epsilon$, the dispersion again gets 
weaker which implies stronger gradients and thus more rapid 
oscillations. In Fig.~\ref{nlsfdseche01s09} we show the solution of 
the semiclassical fNLS equation (\ref{fNLSe}) for $s=0.9$, $\epsilon=0.1$ and 
initial data $ \psi^\e(0,x)\equiv u(x) = \text{sech}(x)$, which is very similar to the 
situation with $s=1$. But it can already be seen here that 
the initial hump splits into two smaller ones, in contrast to the 
case $s=1$. The 
computation is carried out with $N=2^{14}$ Fourier modes and 
$N_{t}=10^{4}$ time steps. 

\begin{figure}[htb!]
   \includegraphics[width=0.49\textwidth]{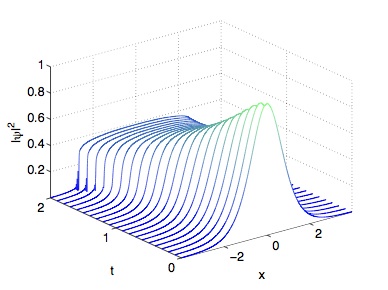}
   \includegraphics[width=0.49\textwidth]{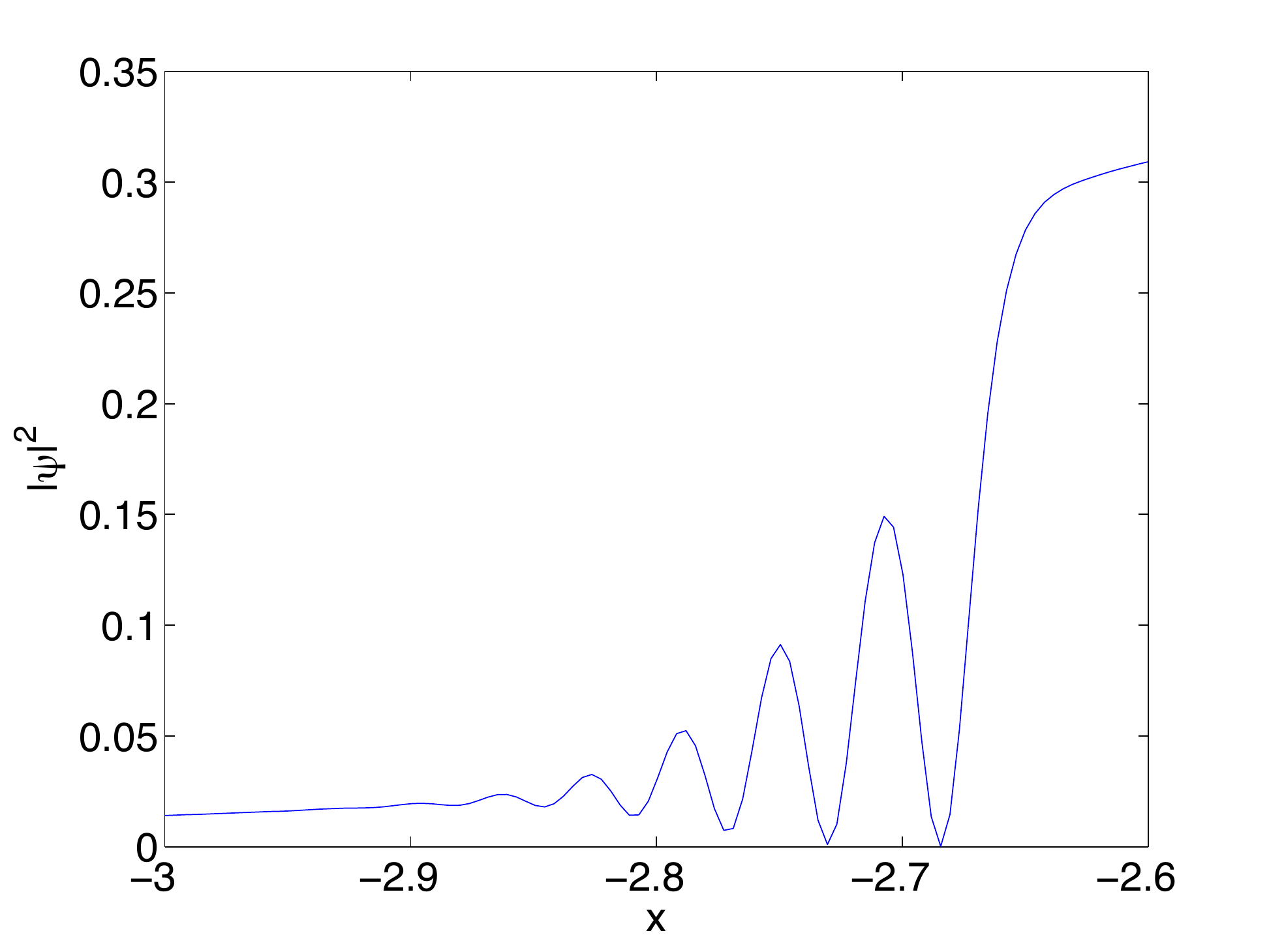}
  \caption{Solution of the  semiclassical, defocusing fNLS equation (\ref{fNLSe})
for the mass subcritical case $s=0.9$, $\epsilon=0.01$, and with initial data 
 $\psi_{0}=\mbox{sech}(x)$; on the right a close-up of the oscillatory 
 zone at the last shown time.}
 \label{nlsfdseche001s09}
\end{figure}

For even smaller $\epsilon$, there are much more oscillations in an 
otherwise identical setting as in Fig.~\ref{nlsfdseche01s09}, as can 
be seen in Fig.~\ref{nlsfdseche001s09}.

\begin{figure}[htb!]
   \includegraphics[width=0.49\textwidth]{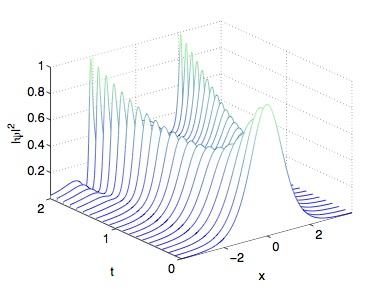}
   \includegraphics[width=0.49\textwidth]{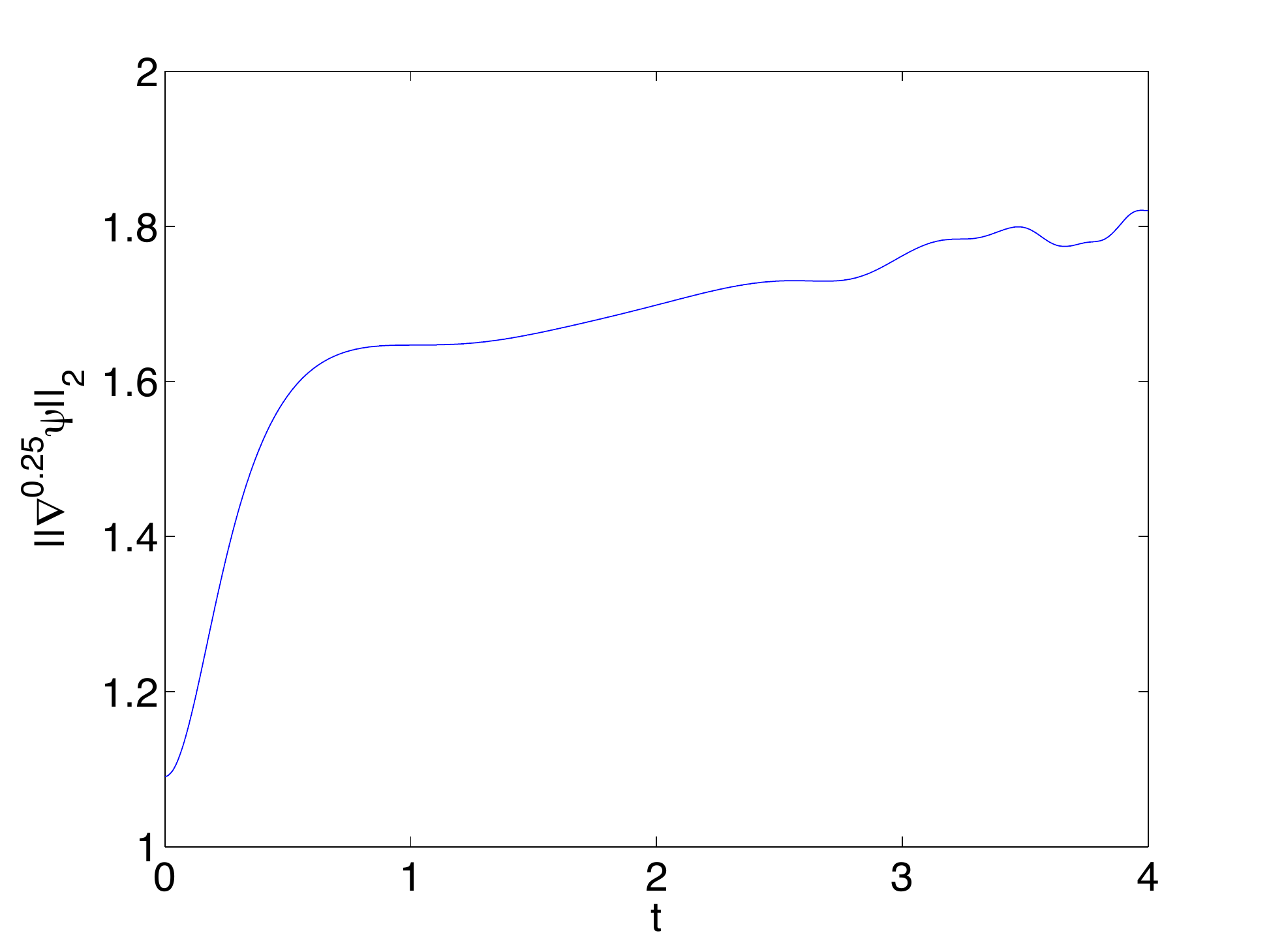}
  \caption{Left: Solution of the semiclassical, defocusing fNLS equation (\ref{fNLSe})
for the energy critical case $s=0.25$, $\epsilon=0.1$, and with initial data 
 $\psi_{0}=\mbox{sech}(x)$. Right: The invariant Sobolev norm 
 (\ref{dotH}); the energy of the initial data implies 
 $\sqrt{2E/\epsilon^{2s}}\approx 2.325$.}
 \label{nlsfdseche01s025}
\end{figure}

\begin{figure}[htb!]
   \includegraphics[width=0.49\textwidth]{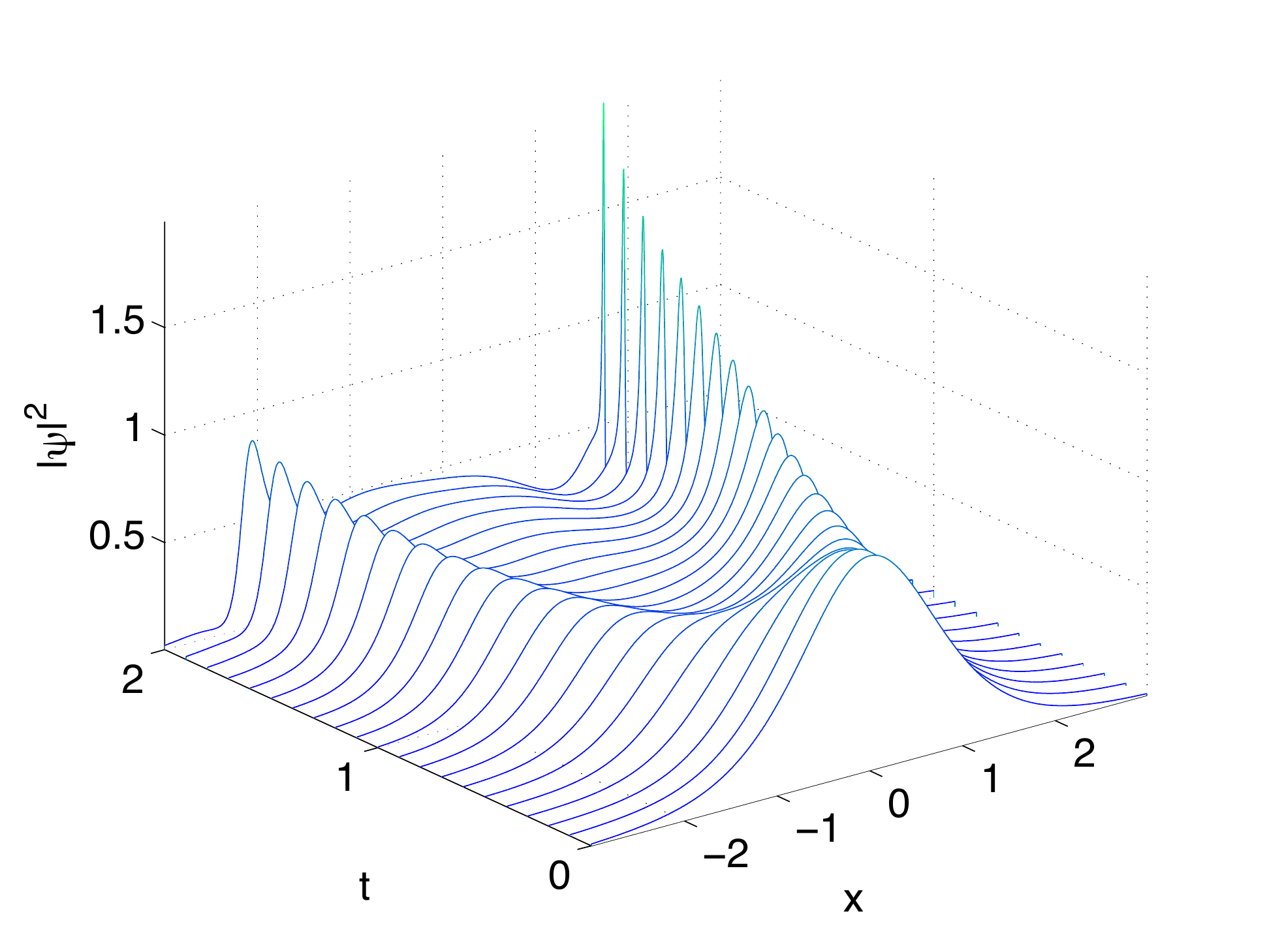}
   \includegraphics[width=0.49\textwidth]{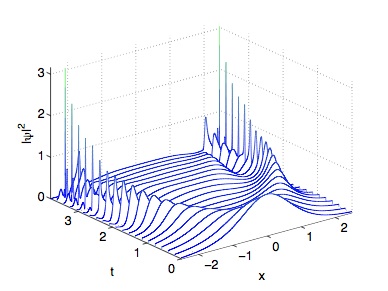}
  \caption{Solution of the semiclassical, defocusing fNLS equation (\ref{fNLSe})
for  $\epsilon=0.1$: on the left for the energy critical case 
$s=0.25$ and initial data 
 $\psi_{0}=e^{ix}\mbox{sech}(x)$, on the right 
 the energy supercritical case $s=0.2$ with initial data 
 $\psi_{0}=\mbox{sech}(x)$.}
 \label{nlsfdseche01s02}
\end{figure}

Reducing $s$ has a similar effect as can be recognized in 
Fig.~\ref{nlsfdseche01s025} where we show the solution for the initial 
data $ \psi^\e(0,x)\equiv u(x) = \text{sech}(x)$ in the energy critical case $s=0.25$. An additional effect of 
the smaller dispersion is that (as noted above) the 
initial hump splits into two humps, which are now well defined (at least before the formation of the first dispersive shock). 
Later in time, there 
appears to be a focusing effect for these two humps, as they get compressed and increase in height. If the code 
is run for longer times, this phenomenon continues and the code 
finally runs out of resolution. We also show the scaling invariant 
$H^\sigma$ norm in the same figure. If we consider instead of the initial data 
 $\psi_{0}=\mbox{sech}(x)$ the same data multiplied by a factor 
 $e^{ix}$, one obtains the solution in Fig.~\ref{nlsfdseche01s02} on 
 the left. The effect already displayed in 
 Fig.~\ref{nlsffracsol_09_eix}, can be also recognized in 
 Fig.~\ref{nlsfdseche01s02}: The hump on the right propagates faster 
 to the right and becomes much earlier `focused' than the one on the 
 left.

To address the question whether the focusing of these humps could 
lead eventually to the formation of a singularity as for solutions of the 
semiclassical system (note that existence of global regular solutions 
is not proven for energy supercritical fNLS), we consider 
the energy supercritical case $s=0.2$ in more detail. The code is run with $N=2^{17}$ for 
$x\in7[\pi,\pi]$ and 
$N_{t}=50000$ time steps for $t\in[0,3.8]$. The solution can be seen 
in Fig.~\ref{nlsfdseche01s02}. The extreme compression of the humps 
is clearly visible.

\begin{figure}[htb!]
   \includegraphics[width=0.49\textwidth]{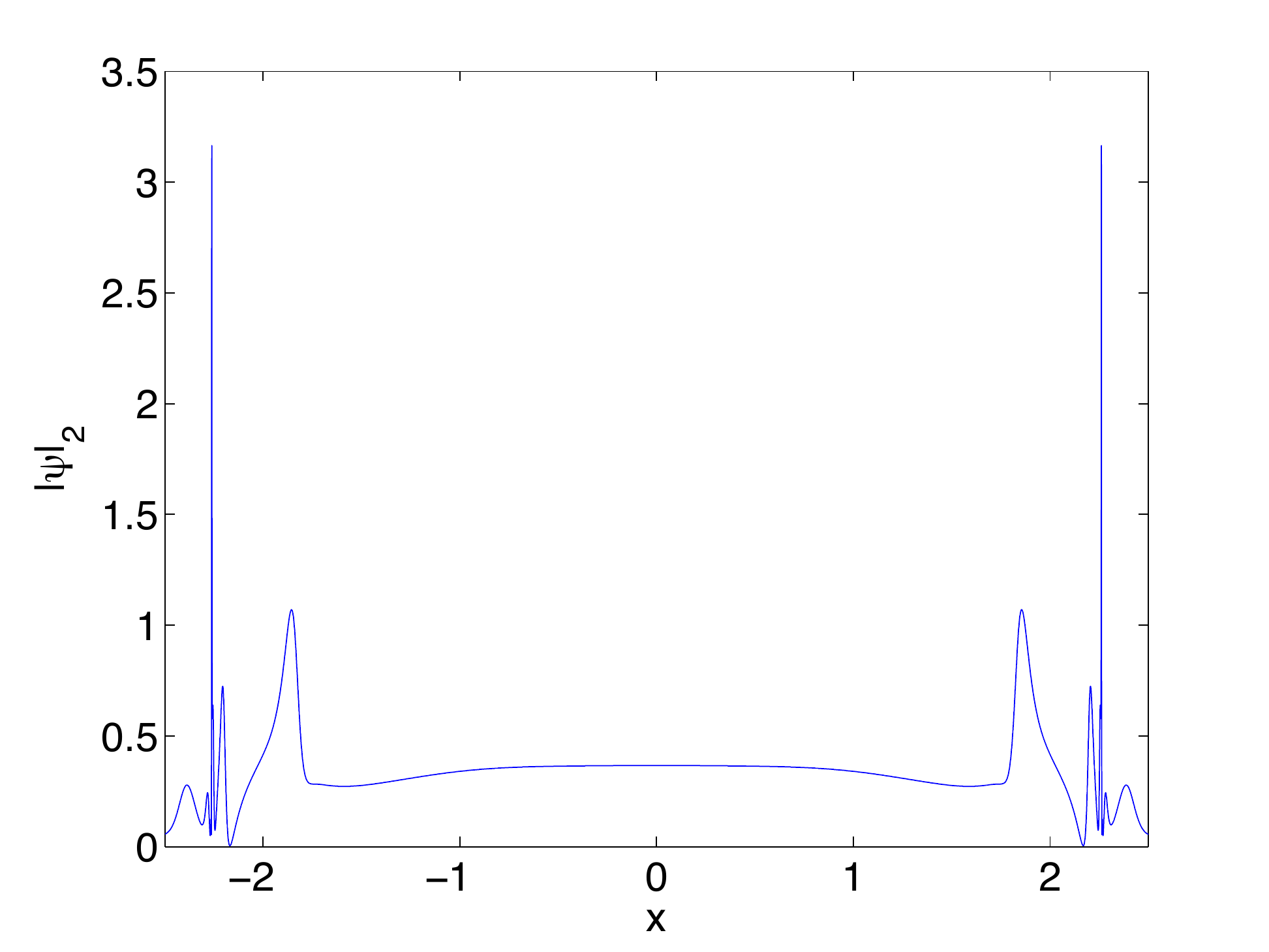}
   \includegraphics[width=0.49\textwidth]{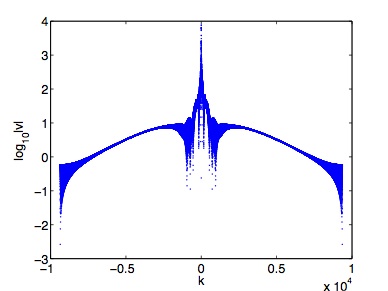}
  \caption{Solution of the semiclassical, defocusing fNLS equation (\ref{fNLSe})
for $s=0.2$, $\epsilon=0.1$ and the initial data 
 $\psi_{0}=\mbox{sech}(x)$ at $t=3.7368$ on the left, and the corresponding 
 Fourier coefficients on the right.}
 \label{nlsfdseche01s02fourier}
\end{figure}

The code is stopped at $t=3.7368$ since the distance of the nearest 
singularity in the complex plane to the real axis as determined by 
fitting the Fourier coefficients to the asymptotic formula 
(\ref{fourasymp}) is smaller than 
the smallest resolved distance in physical space. But as can be seen 
from Fig.~\ref{nlsfdseche01s02fourier}, where we also show the 
solution at the last recorded time, the Fourier coefficients indicate 
that the code ran out of resolution before. In fact, the modulus of 
the Fourier coefficients decreases only to the order of $10^{-1}$ at 
$t=3.61$, whereas it reached $10^{-6}$ for $t=3.41$ (this implies 
that the numerical results in this case should be ignored for 
$t>3.5$). Thus the code 
runs out of resolution well before a potential singularity hits the 
real axis. Rerunning the code with higher resolution produces the 
same phenomena, just at slightly later times. 
This indicates that the solutions indeed stays regular in this case, but is nevertheless very 
different from the focusing fNLS equation studied in the previous subsection. It is also different 
from the well-known cusps found in the case of semiclassical defocusing NLS, as can be  
identified using the techniques given in 
\cite{KR2013b}. We finally note that the same qualitative behavior is also observed 
for different choices of localized initial data.

\begin{figure}[htb!]
   \includegraphics[width=0.32\textwidth]{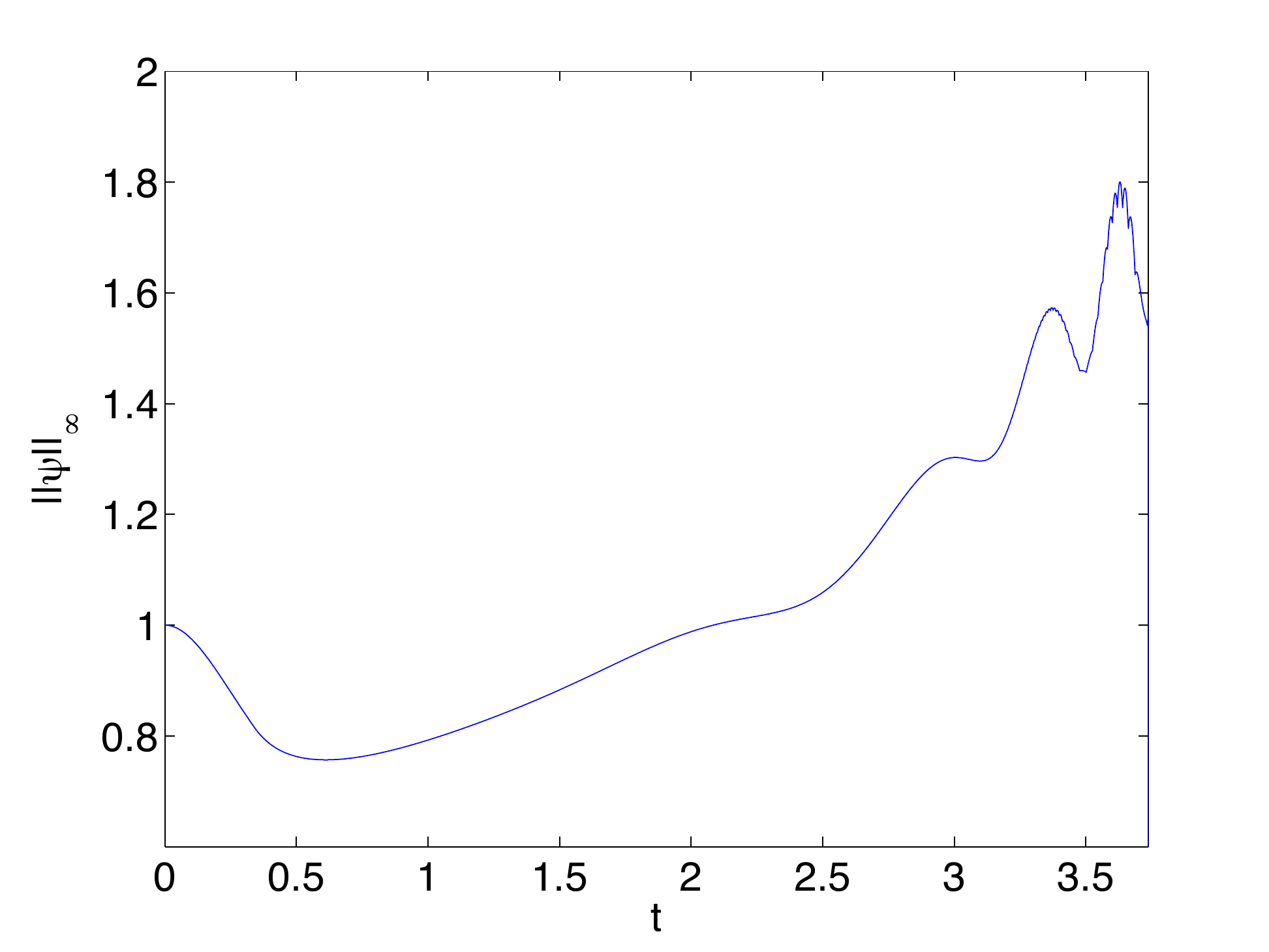}
   \includegraphics[width=0.32\textwidth]{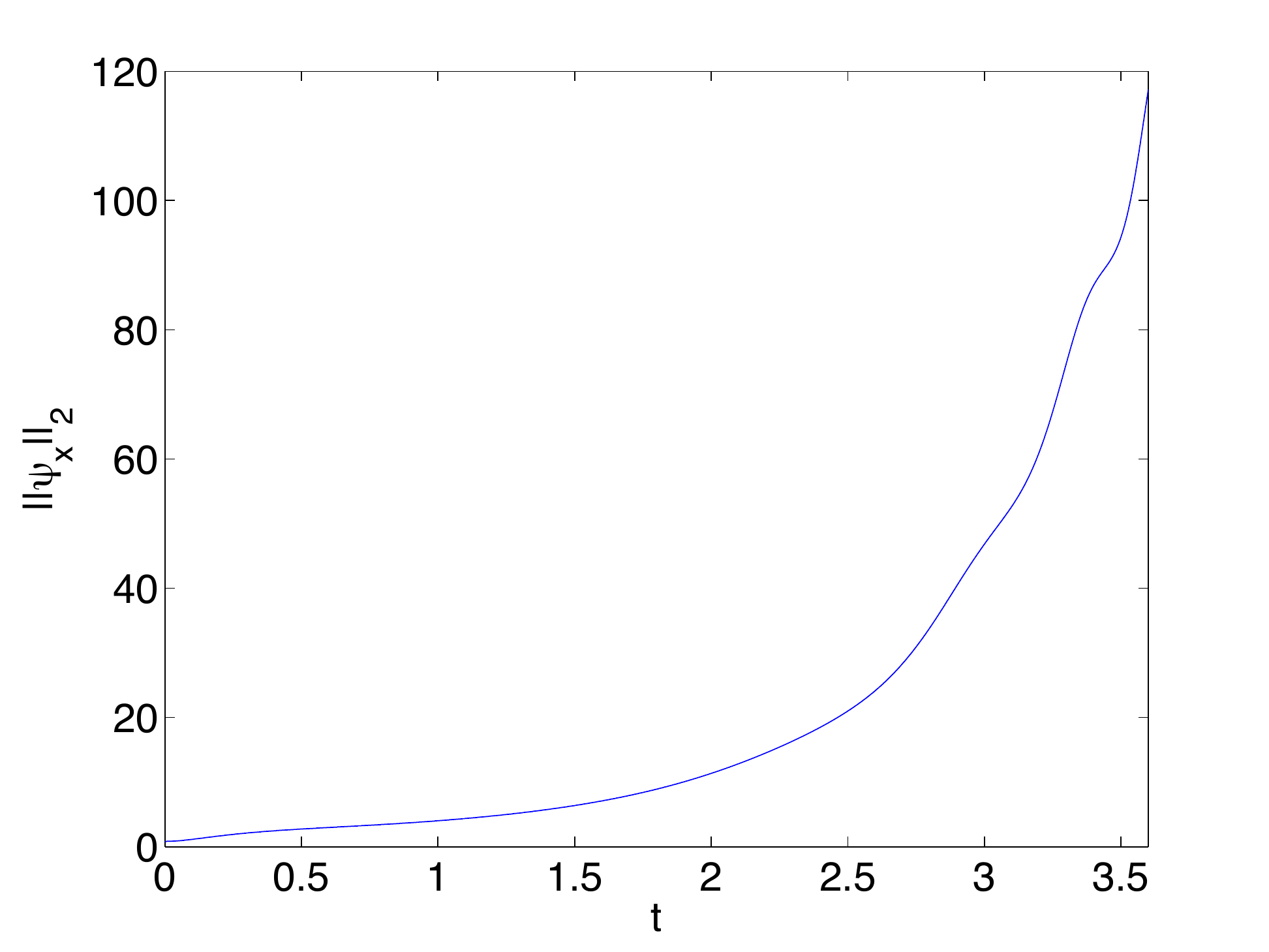}
   \includegraphics[width=0.32\textwidth]{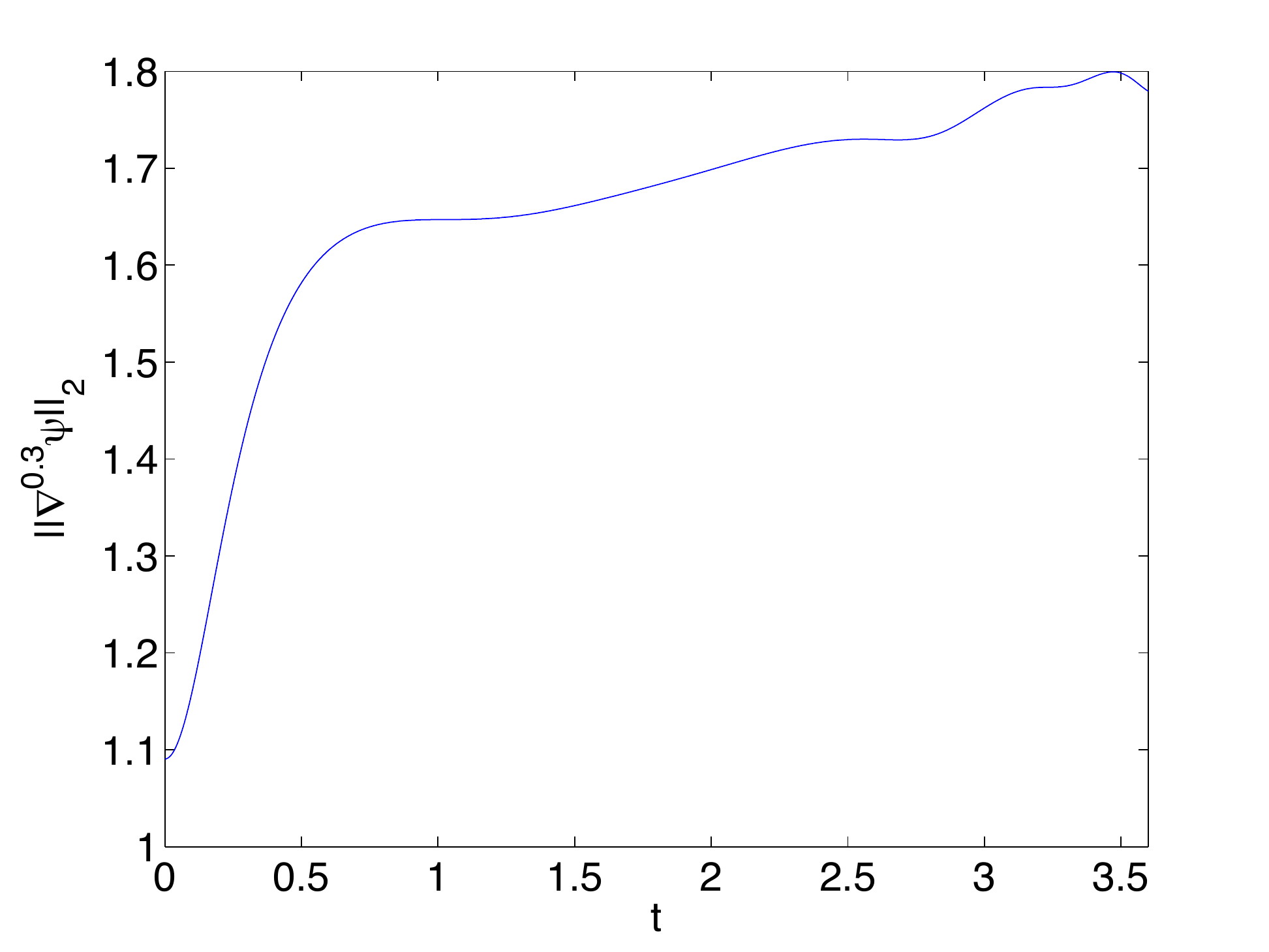}
  \caption{Time-dependence of various norms of the solution of the semiclassical, defocusing fNLS equation (\ref{fNLSe})
for $s=0.2$, $\epsilon=0.1$ and the initial data 
 $\psi_{0}=\mbox{sech}(x)$:  on the left the 
 $L^\infty$ norm, in the middle the $L^2$ norm of the gradient, 
 and on the right the invariant norm (\ref{dotH}).}
 \label{nlsfdseche01s02norm}
\end{figure}

Various norms of this solution are shown in 
Fig.~\ref{nlsfdseche01s02norm}, where one has to bear in mind that 
there is a lack of resolution for 
the last time steps. It can be seen that the $L^\infty$ norm continues 
to grow (the shown oscillations might be spurious and due to finite 
resolution in physical space), and that there is a strong growth in the $L^2$ norm of 
the gradient of the solution. Whereas the strong gradients in the 
solution are reflected by the growth of the latter norm, there 
is no indication of it blowing up at a time close to the last 
computed time. Also the $\dot H^\sigma$ norm invariant under the 
rescaling (\ref{scaling}) grows only moderately. This is in 
accordance with the conclusion above that the solution stays regular. 
Note that the strong compression of the humps visible in 
Fig.~\ref{nlsfdseche01s02} does not allow us to reach the same 
asymptotic regime for the long time behavior of the solution as in 
\cite{CSS} for a higher dimensional NLS equation. The reason for this 
is simply that the solutions decay much slower in $|x|$ in one 
spatial dimension, especially in the presence of fractional 
derivatives. 


\section{Conclusion}\label{sec:conclusio}
In this paper we have presented a comprehensive numerical study of 
 issues appearing in the context of fractional NLS equations 
in one spatial dimension. We have fixed the nonlinearity to be cubic 
and have varied the order $s$ of the fractional derivatives. This allowed 
to explore the competing effects of nonlinearity and dispersion 
in NLS systems. 

Concretely we were able to numerically construct ground state 
solutions of the 
focusing fNLS equation and study their stability in certain regimes. As expected, the 
ground states are stable for $s>\frac{1}{2}$ (the mass subcritical case). Perturbations of 
such states result in a solution which numerically displays damped 
oscillations around a final ground state of the same mass as the 
initial data. 
Moreover, we also find that approximate solitary wave solutions are possible despite the 
non-locality of the dispersion within our model. For smaller values of $s\le \frac{1}{2}$, the 
ground state is unstable both against being radiated away towards 
infinity and finite-time blow-up. 
Concerning the latter, we numerically studied the appearance of blow-up in mass and energy supercritical regimes.
In the mass critical regime, we give numerical evidence for a self-similar blow-up with a profile given 
by the fractional ground state $Q(x)$ and a rate close to the 
well-known log-log regime (as proved for mass critical NLS with $s=1$ and $\sigma =\frac{2}{d}$). 

We also studied the long time behavior of solutions in the energy 
critical and supercritical regime. To this end we introduced a 
semiclassical parameter $\e>0$ in the equation (\ref{fNLSe}) and studied the 
corresponding asymptotic regime for $\e \ll 1$. We found that 
solutions of the defocusing fNLS 
equation appear to stay regular in time but exhibit a surprising oscillatory structure which is much more extreme than 
for the associated NLS ($s=1$) case.
In the focusing case, blow-up was found for values of $s\le \frac{1}{2}$, but there seems to be no qualitative difference between 
mass-supercritical and energy-supercritical blow-up.

%
%
%

\end{document}